\documentclass[11pt]{amsart}
\usepackage{amsmath,amsthm,amssymb}
\usepackage{mathrsfs}
\usepackage[titletoc]{appendix}
\usepackage{bbm}
\usepackage{float}
\usepackage{tikz}
\begin{document}
\input amssym.def

\setcounter{equation}{0}
\newcommand{\wt}{\mbox{wt}}
\newcommand{\spa}{\mbox{span}}
\newcommand{\Res}{\mbox{Res}}
\newcommand{\End}{\mbox{End}}
\newcommand{\Ind}{\mbox{Ind}}
\newcommand{\Hom}{\mbox{Hom}}
\newcommand{\Mod}{\mbox{Mod}}
\newcommand{\m}{\mbox{mod}\ }
\newcommand{\ch}{\mbox{ch}}
\renewcommand{\theequation}{\thesection.\arabic{equation}}
\numberwithin{equation}{section}

\def \Aut{{\rm Aut}}
\def \Ad{{\rm Ad}}
\def \r{\rho}
\def \Z{\Bbb Z}
\def \M{\Bbb M}
\def \C{\Bbb C}
\def \R{\Bbb R}
\def \P{\Bbb P}
\def \Q{\Bbb Q}
\def \N{\Bbb N}
\def \ann{{\rm Ann}}
\def \<{\langle}
\def \Om{\Omega}
\def \M{{\cal M}}
\def \1t{\frac{1}{T}}
\def \>{\rangle}
\def \t{\tau }
\def \a{\alpha }
\def \d{\delta}
\def \e{\epsilon }
\def \l{\lambda }
\def \L{\Lambda }
\def \g{\gamma}
\def \b{\beta }
\def \ro {\rho}
\def \o {\omega}
\def \I {\mathcal{I}}
\def \A {\mathcal{A}}
\def \B {\mathcal {B}}
\def \Cc {\mathcal {C}}
\def \H {\mathcal{H}}
\def \M {\mathcal{M}}
\def \V {\mathcal{V}}
\def \Y {\mathcal{Y}}
\def \Mob {\hbox{M\"{o}b}}
\def \res {|}
\def \om{\omega }
\def \o{\omega }
\def \cg{\chi_g}
\def \ag{\alpha_g}
\def \ah{\alpha_h}
\def \ph{\psi_h}
\def \nor{\vartriangleleft}
\def \Diff {\hbox {Diff}}
\def \col {\hbox{col}}
\def \Hom {\hbox{Hom}}
\def \la {\langle}
\def \ra {\rangle}
\def \voa{vertex operator algebra }
\def \voas{vertex operator algebras  }
\def \vosa{vertex operator subalgebra }
\def \v{vertex operator algebra\ }
\def \1{{\bf 1}}
\def \be{\begin{equation}\label}
\def \ee{\end{equation}}
\def \qed{\mbox{ $\square$}}
\def \pf {\noindent {\bf Proof:} \,}
\def \bl{\begin{lem}\label}
\def \el{\end{lem}}
\def \ba{\begin{array}}
\def \ea{\end{array}}
\def \bt{\begin{thm}\label}
\def \et{\end{thm}}
\def \ch{{\rm ch}}
\def \br{\begin{rem}\label}
\def \er{\end{rem}}
\def \ed{\end{de}}
\def \bp{\begin{prop}\label}
\def \ep{\end{prop}}

\newtheorem{thm}{Theorem}[section]
\newtheorem{prop}[thm]{Proposition}
\newtheorem{lem}[thm]{Lemma}
\newtheorem{cor}[thm]{Corollary}
\newtheorem{rmk}[thm]{Remark}
\newtheorem*{CPM}{Mirror Extension Conjecture}
\newtheorem{definition}[thm]{Definition}
\newtheorem{conj}[thm]{Conjecture}
\newtheorem{ex}[thm]{Example}
%%%%%%%%%%%%%%%%%%%%%%%%%%%%%%%%%%%%%%%%%%%%%%%%%%%%%

\title
{Mirror Extensions of Vertex Operator Algebras}

\author{Chongying Dong}
\address{Department of Mathematics, Sichuan University,
 Chengdu, China \& Department of Mathematics, University of California, Santa Cruz 98064, USA}
\email{dong@ucsc.edu}
 \thanks{The first author was partially supported by NSF grants}
\author{Xiangyu Jiao}
\address{Department of Mathematics, University of California, Santa Cruz 98064, USA}
\email{xjiao@ucsc.edu}
\author{Feng Xu}
\address{Department of Mathematics, University of California, Riverside, CA 92521, USA}
\email{xufeng@math.ucr.edu}
 \thanks{The third author was partially supported by NSF grants}

\begin{abstract}
The mirror extensions for vertex operator algebras are studied. Two explicit examples which
are not simple current extensions of some affine vertex operator algebras of type $A$ are given.
\end{abstract}

\maketitle
%---------------------------------------------------------------------------
\section{Introduction}
Mirror extensions, in the title of this paper, refer to a general Theorem 3.8 in \cite{Xu-m} which produces
completely rational conformal nets from given ones. Based on the close relations between conformal
nets and vertex operator algebras, we make the following conjecture
which is the vertex operator algebra version of Theorem 3.8 in \cite{Xu-m}:

\begin{CPM}\label{conj}
Let $V$ be a rational and $C_2$-cofinite \voa and $U$  a rational and $C_2$-cofinite \vosa of $V.$
Denote $U^c$  the commutant \voa of $U$ in $V.$
 Assume that $(U^c)^c=U,$ and
 $$V=U\otimes U^c\bigoplus (\oplus_{i=1}^n U_i\otimes U^c_i),$$
  where $U_i$'s and $U^c_i$'s are irreducible modules for $U$ and $U^c$ respectively.
 Then if
 $$ U^e=U\bigoplus (\oplus_{i=1}^n m_iU_{i})$$
 is a rational \voa where $m_i\geq 0,$ so is
  $$(U^c)^e=U^c\bigoplus (\oplus_{i=1}^n m_iU^c_{i}).$$
\end{CPM}

 The mirror extensions of conformal nets associated to affine Kac-Moody
 algebras of type $A$ have been studied extensively in \cite{Xu-m}.
 It is conjectured on P.846 of \cite{Xu-m} that there should be rational
vertex operator algebras corresponding to the class of completely rational conformal nets constructed
 in $\S4.3$ of \cite{Xu-m}, and this conjecture which is a special case of Mirror Extension Conjecture
 is the  motivation of our paper.

A proof of these conjectures seems to be out of reach at present. Instead we focus
 %translation of Theorem 3.8 in \cite{Xu-m} from conformal net language to vertex operator algebra
on two interesting examples of mirror extensions of
\voas(cf. P.836 of \cite{Xu-m}) in this paper. The first example is based on conformal inclusions
$SU(2)_{10} \subset Spin(5)_1$ and $SU(2)_{10} \times SU(10)_2\subset SU(20)_1$. The spectrum of
$SU(2)_{10} \subset Spin(5)_1$ is $H_0+H_6$, and $(6, \Lambda_3+\Lambda_7)$ appears in the spectrum
of  $SU(2)_{10} \times SU(10)_2\subset SU(20)_1.$ Here we use $\Lambda_i$ to denote the fundamental
weights of $SU(n),$ and $0$ (or $\Lambda_0$) to denote the trivial representation of $SU(n)$ and we specialize the case
$SU(2)$ by using $i$ to denote the highest weight of the representation of $SU(2)$.
Theorem 3.8 in \cite{Xu-m} implies that there exists a completely rational net containing $\mathcal{A}_{SU(10)_2}$
with spectrum $H_{0}+H_{\Lambda_3+\Lambda_7}$.

From \voa point of view, this suggests that there should be a \voa structure on
$L_{\frak{sl}(10)}(2,0)+L_{\frak{sl}(10)}(2,\Lambda_3+\Lambda_7)$ where $L_{\frak{g}}(k,\lambda)$ is the
highest weight irreducible module for the affine Lie algebra $\hat{\frak g}$ of level $k$ associated to the weight
$\lambda$ of $\frak{g}.$   Notice that the lowest weight
$h_{\Lambda_3+\Lambda_7}$ of $L_{\frak{sl}(10)}(2,\Lambda_3+\Lambda_7)$ is $2$ and
$L_{\frak{sl}(10)}(2,\Lambda_3+\Lambda_7)$ is not a simple current, such a \voa has not been obtained
from the affine \voa in the literature. We need to find
intertwining operators associated to vectors in $L_{\frak{sl}(10)}(2,\Lambda_3+\Lambda_7)$
 verifying the crucial locality condition of vertex operator algebra. The correlation functions of such intertwining operators  are
  solutions of KZ equation, and locality means that such functions are symmetric rational functions. Thus in
  this case we need to find symmetric rational solutions to KZ equation for $SU(10)_2$.

There are contour integral representations of solutions  to KZ equation for $SU(n)_k$ (cf. \cite{FEK} and references therein) for generic level $k\notin \Q.$ It is not clear how to find general solutions and pick out a particular rational
solution in our case with $n=10$ and $k=2.$ Instead we take a different approach, which contains the key
idea of this paper. First we note that from the conformal inclusion $SU(2)_{10}\subset Spin(5)$, the
primary fields in $L_{\frak{sl}(2)}(10,6)$ in this inclusion produce correlator X which is a symmetric rational solution of KZ equation for
the affine \voa $L_{\frak{sl}(2)}(10,0)$. This is equivalent to say that these solutions are invariant under braiding
operator B. By abusing of notation we write $BX=X.$ Note that $(6,\Lambda_3+\Lambda_7)$ appears in the
spectrum of $SU(2)_{10}\times SU(10)_2\subset SU(20)_1$, the vertex operator associated to the highest
weight vector of $(6,\Lambda_3+\Lambda_7)$ will give us a symmetric rational function. Due to a crucial
non-degenerate property in Corollary \ref{nondeg}, this implies $B'\dot{B}=Id$, where $B'$ is similar to $B$
by conjugation of invertible diagonal matrix. Due to the  crossing symmetry property of $B$ for $SU(2)$ in Lemma \ref{cross sym}
from $BX=X$ we conclude there must be $\dot{X}$ which verifies KZ equation for $L_{\frak{sl}(10)}(2,0)$
such that $\dot{B}\dot{X}=\dot{X}$,
and it follows that $\dot{X}$ is a symmetric rational function.

From $\dot{X}$ we can easily define a \voa structure on
$$L_{\frak{sl}(10)}(2,0)\oplus L_{\frak{sl}(10)}(2,\Lambda_3+\Lambda_7),$$ and we derive our main
result Theorem \ref{main}. The \voa in Theorem \ref{main} is an example of mirror extension, which is constructed
by an idea very different from what is previously known. We are informed recently that the \voa
$L_{\frak{sl}(10)}(2,0)\oplus L_{\frak{sl}(10)}(2,\Lambda_3+\Lambda_7)$ can also be realized a coset
construction  in a holomorphic vertex operator algebra with central charge $24$ \cite{L}.

The second example is based on the conformal inclusion $SU(2)_{28}\subset (G_2)_1$ (see \cite{CIZ, GNO}) and
the level-rank duality $SU(2)_{28}\times SU(28)_2\subset SU(56)_1.$ Similar to the first example,
$$L_{\frak{sl}(28)}(2,0)\oplus L_{\frak{sl}(28)}(2,\Lambda_5+\Lambda_{23})\oplus L_{\frak{sl}(28)}(2,\Lambda_9+\Lambda_{19})
\oplus L_{\frak{sl}(28)}(2,2\Lambda_{14})$$
is a vertex operator algebra which is a mirror extension corresponding to the vertex operator algebra
$$L_{G_2}(1,0)=L_{\frak{sl}(2)}(28,0)\oplus L_{\frak{sl}(2)}(28,10)\oplus L_{\frak{sl}(2)}(28,18)\oplus L_{\frak{sl}(2)}(28,28).$$

Although these two examples of mirror extensions which are not simple current extensions are totally new in the theory of vertex operator algebra, the mirror extension, in fact, is a general phenomenon. Many well known vertex operator algebras in the literature can also be regraded as mirror extensions. We give two easy examples here. The first example  comes from the well known GKO-construction \cite{GKO}:
 \begin{equation*}
 \begin{aligned}
 L_{\frak{sl}(2)}(1,0) \otimes L_{\frak{sl}(2)}(3,0)&=L(\frac{4}{5},0)\otimes L_{\frak{sl}(2)}(4,0)
 \oplus L(\frac{4}{5},\frac{2}{3})\otimes L_{\frak{sl}(2)}(4,2)\\
&\oplus L(\frac{4}{5},3)\otimes L_{\frak{sl}(2)}(4,4),
 \end{aligned}
 \end{equation*}
 where $L(\frac{4}{5},h)$ is the lowest weight irreducible module for the
 Virasoro algebra with central charge $\frac{4}{5}$ and lowest weight $h.$
The \voa structure on $L(\frac{4}{5},0)\oplus L(\frac{4}{5},3)$ is well known now (see \cite{KMY}),
 which is a simple current extension of $L(\frac{4}{5},0).$ The vertex operator algebra $L_{\frak{sl}(2)}(4,0)\oplus L_{\frak{sl}(2)}(4,4)$ (see \cite{MS, Li2}) is a mirror extension.

The other is the $3A$ algebra $U$ \cite{LYY} which has a decomposition:
\begin{equation*}
\begin{aligned}
U\cong &L(\frac{4}{5},0)\otimes L(\frac{6}{7},0)\oplus L(\frac{4}{5},3)\otimes  L(\frac{6}{7},5)
\oplus L(\frac{4}{5},\frac{2}{3})\otimes  L(\frac{6}{7},\frac{4}{3})\\
&\oplus L(\frac{4}{5},\frac{13}{8})\otimes  L(\frac{6}{7},\frac{3}{8})
\oplus L(\frac{4}{5},\frac{1}{8})\otimes  L(\frac{6}{7},\frac{23}{8}).
\end{aligned}
\end{equation*}
Again the vertex operator algebra $ L(\frac{6}{7},0)\oplus  L(\frac{6}{7},5)$ \cite{LY} is
a mirror extension corresponding to the vertex operator algebra  $L(\frac{4}{5},0)\oplus L(\frac{4}{5},3).$

Besides what is already described above, we have included a preliminary section
$\S$\ref{prelim} on affine vertes operator algebras, KZ equation, primary fields, conformal nets and induction, and we prove the crucial non-degeneracy condition in Corollary \ref{nondeg}.
The first nontrivial example of mirror
extensions  is presented in $\S$\ref{key}. The uniqueness of this \voa structure
is obtained in $\S$\ref{uni}. In $\S$\ref{28} we construct another example by using similar method.
In $\S$\ref{comments} we discuss problems about general case. The last section is the appendix which is devoted to proving
the non-degeneracy property given in Corollary \ref{nondeg} using \voa language.
%----------------------------------------------------------------------------------
%----------------------------------------------------------------------------------

\section{Preliminaries}\label{prelim}

\subsection{Preliminaries on sectors}

Given an infinite factor $M$, the {\it sectors of $M$}  are given by
$$\text{Sect}(M) = \text{End}(M)/\text{Inn}(M),$$
namely $\text{Sect}(M)$ is the quotient of the semigroup of the
endomorphisms of $M$ modulo the equivalence relation: $\rho,\rho'\in
\text{End}(M),\, \rho\thicksim\rho'$ iff there is a unitary $u\in M$
such that $\rho'(x)=u\rho(x)u^*$ for all $x\in M$.

$\text{Sect}(M)$ is a $^*$-semiring (there are an addition, a
product and an involution $\rho\rightarrow \bar\rho$) equivalent to
the Connes correspondences (bimodules) on $M$ up to unitary
equivalence. If $\ro$ is an element of $\text{End}(M)$ we shall
denote by $[\ro]$ its class in $\text{Sect}(M)$. We define
$\text{Hom}(\ro,\ro')$ between the objects $\ro,\ro'\in \End(M)$ by
\[
\text{Hom}(\ro,\ro')\equiv\{a\in M: a\ro(x)=\ro'(x)a \ \forall x\in M\}.
\]
We use $\langle  \lambda , \mu \rangle$ to denote the dimension of
$\text{\rm Hom}(\lambda , \mu )$; it can be $\infty$, but it is
finite if $\l,\mu$ have finite index. See \cite{J1} for the
definition of index for type $II_1$ case which initiated the subject
and  \cite{PP} for  the definition of index in general. Also see
\S2.3 \cite{KLX} for expositions. $\langle  \lambda , \mu \rangle$
depends only on $[\lambda ]$ and $[\mu ]$. Moreover we have if $\nu$
has finite index, then $\langle \nu \lambda , \mu \rangle = \langle
\lambda , \bar \nu \mu \rangle $, $\langle \lambda\nu , \mu \rangle
= \langle \lambda , \mu \bar \nu \rangle $ which follows from
Frobenius duality. $\mu $ is a subsector of $\lambda $ if there is
an isometry $v\in M$ such that $\mu(x)= v^* \lambda(x)v, \forall
x\in M.$ We will also use the following notation: if $\mu $ is a
subsector of $\lambda $, we will write as $\mu \prec \lambda $  or
$\lambda \succ \mu $.  A sector is said to be irreducible if it has
only one subsector.

\subsection{Local nets}\label{ci}
%\subsubsection{Preliminaries on conformal nets}
By an interval of the circle we mean an open connected non-empty
subset $I$ of $S^1$ such that the interior of its complement $I'$ is
not empty. We denote by $\I$ the family of all intervals of $S^1$.
A {\it net} $\A$ of von Neumann algebras on $S^1$ is a map
\[
I\in\I\to\A(I)\subset B(\H)
\]
from $\I$ to von Neumann algebras on a fixed separable Hilbert space
$\H$ that satisfies:
\begin{itemize}
\item[{\bf A.}] {\it Isotony}. If $I_{1}\subset I_{2}$ belong to
$\I$, then
\begin{equation*}
 \A(I_{1})\subset\A(I_{2}).
\end{equation*}
\end{itemize}
If $E\subset S^1$ is any region, we shall put
$\A(E)\equiv\bigvee_{E\supset I\in\I}\A(I)$ with $\A(E)=\mathbb C$
if $E$ has empty interior (the symbol $\vee$ denotes the von Neumann
algebra generated).

The net $\A$ is called {\it local} if it satisfies:
\begin{itemize}
\item[{\bf B.}] {\it Locality}. If $I_{1},I_{2}\in\I$ and $I_1\cap
I_2=\emptyset$ then
\begin{equation*}
 [\A(I_{1}),\A(I_{2})]=\{0\},
 \end{equation*}
where brackets denote the commutator.
\end{itemize}
The net $\A$ is called {\it M\"{o}bius covariant} if in addition
satisfying the following properties {\bf C,D,E,F}:
\begin{itemize}
\item[{\bf C.}] {\it M\"{o}bius covariance}.
There exists a non-trivial strongly continuous unitary
representation $U$ of the M\"{o}bius group $\Mob$ (isomorphic to
$PSU(1,1)$) on $\H$ such that
\begin{equation*}
 U(g)\A(I) U(g)^*\ =\ \A(gI),\quad g\in \Mob,\ I\in\I.
\end{equation*}
\item[{\bf D.}] {\it Positivity of the energy}.
The generator of the one-parameter rotation subgroup of $U$
(conformal Hamiltonian), denoted by $L_0$ in the following,  is
positive.
\item[{\bf E.}] {\it Existence of the vacuum}.  There exists a unit
$U$-invariant vector $\Omega\in\H$ (vacuum vector), and $\Omega$ is
cyclic for the von Neumann algebra $\bigvee_{I\in\I}\A(I)$.
\end{itemize}
By the Reeh-Schlieder theorem $\Omega$ is cyclic and separating for
every fixed $\A(I)$. The modular objects associated with
$(\A(I),\Omega)$ have a geometric meaning
\[
\Delta^{it}_I = U(\Lambda_I(2\pi t)),\qquad J_I = U(r_I)\ .
\]
Here $\Lambda_I$ is a canonical one-parameter subgroup of $\Mob$ and
$U(r_I)$ is an antiunitary acting geometrically on $\A$ as a
reflection $r_I$ on $S^1$.

This implies {\em Haag duality}:
\[
\A(I)'=\A(I'),\quad I\in\I\ ,
\]
where $I'$ is the interior of $S^1\setminus I$.

\begin{itemize}
\item[{\bf F.}] {\it Irreducibility}. $\bigvee_{I\in\I}\A(I)=B(\H)$.
Indeed $\A$ is irreducible iff $\Om$ is the unique $U$-invariant
vector (up to scalar multiples). Also  $\A$ is irreducible iff the
local von Neumann algebras $\A(I)$ are factors. In this case they
are either ${\mathbb C}$ or III$_1$-factors with separable predual
in Connes classification of type III factors.
\end{itemize}
By a {\it conformal net} (or diffeomorphism covariant net) $\A$ we
shall mean a M\"{o}bius covariant net such that the following holds:
\begin{itemize}
\item[{\bf G.}] {\it Conformal covariance}. There exists a projective
unitary representation $U$ of $\Diff(S^1)$ on $\H$ extending the
unitary representation of $\Mob$ such that for all $I\in\I$ we have
\begin{gather*}
 U(\phi)\A(I) U(\phi)^*\ =\ \A(\phi.I),\quad  \phi\in\Diff(S^1), \\
 U(\phi)xU(\phi)^*\ =\ x,\quad x\in\A(I),\ \phi\in\Diff(I'),
\end{gather*}
\end{itemize}
where $\Diff(S^1)$ denotes the group of smooth, positively oriented
diffeomorphism of $S^1$ and $\Diff(I)$ the subgroup of
diffeomorphisms $g$ such that $\phi(z)=z$ for all $z\in I'$.
\par
A (DHR) representation $\pi$ of $\A$ on a Hilbert space $\H$ is a
map $I\in\I\mapsto  \pi_I$ that associates to each $I$ a normal
representation of $\A(I)$ on $B(\H)$ such that
\[
\pi_{\tilde I}\res\A(I)=\pi_I,\quad I\subset\tilde I, \quad I,\tilde
I\subset\I\ .
\]
$\pi$ is said to be M\"obius (resp. diffeomorphism) covariant if
there is a projective unitary representation $U_{\pi}$ of $\Mob$
(resp. $\Diff(S^1)$) on $\H$ such that
\[
\pi_{gI}(U(g)xU(g)^*) =U_{\pi}(g)\pi_{I}(x)U_{\pi}(g)^*
\]
for all $I\in\I$, $x\in\A(I)$ and $g\in \Mob$ (resp.
$g\in\Diff(S^1)$).

By definition the irreducible conformal net is in fact an
irreducible representation of itself and we will call this
representation the {\it vacuum representation}.\par

Let $G$ be a simply connected  compact Lie group. By Th. 3.2 of
\cite{FG}, the vacuum positive energy representation of the loop
group $LG$ (cf. \cite{PS}) at level $k$ gives rise to an irreducible
conformal net denoted by {\it ${\A}_{G_k}$}. By Th. 3.3 of
\cite{FG}, every irreducible positive energy representation of the
loop group $LG$ at level $k$ gives rise to  an irreducible covariant
representation of ${\A}_{G_k}$. \par Given an interval $I$ and a
representation $\pi$ of $\A$, there is an {\em endomorphism of $\A$
localized in $I$} equivalent to $\pi$; namely $\ro$ is a
representation of $\A$ on the vacuum Hilbert space $\H$, unitarily
equivalent to $\pi$, such that
$\ro_{I'}=\text{id}\restriction\A(I')$. We now define  the
statistics. Given the endomorphism $\ro$ of $\A$ localized in
$I\in\I$, choose an equivalent endomorphism $\ro_0$ localized in an
interval $I_0\in\I$ with $\bar I_0\cap\bar I =\emptyset$ and let $u$
be a local intertwiner in $\Hom(\ro,\ro_0)$ , namely $u\in
\Hom(\ro_{\tilde I},\ro_{0,\tilde I})$ with $I_0$ following clockwise
$I$ inside $\tilde I$ which is an interval containing both $I$ and
$I_0$.

The {\it statistics operator} $\epsilon (\ro,\rho):= u^*\ro(u) =
u^*\ro_{\tilde I}(u) $ belongs to $\Hom(\ro^2_{\tilde I},\ro^2_{\tilde
I})$. We will call $\epsilon (\ro,\rho)$ the positive or right
braiding and $\tilde\epsilon (\ro,\rho):=\epsilon (\ro,\rho)^*$ the
negative or left braiding.
\par
Let $\B$ be a  conformal  net. By a {\it  conformal subnet} (cf.
\cite{Longo4}) we shall mean a map
\[
I\in\I\to\A(I)\subset \B(I)
\]
that associates to each interval $I\in \I$ a von Neumann subalgebra
$\A(I)$ of $\B(I)$, which is isotonic
\[
\A(I_1)\subset \A(I_2), I_1\subset I_2,
\]
and  conformal covariant with respect to the representation $U$,
namely
\[
U(g) \A(I) U(g)^*= \A(g.I)
\] for all $g\in \Diff(S^1)$ and $I\in \I$. Note that by Lemma 13
of \cite{Longo4} for each $I\in \I$ there exists a conditional
expectation $E_I: \B(I)\rightarrow \A(I)$ such that $E_I$ preserves
the vector state given by the vacuum of $\A$.
\begin{definition}\label{ext}
Let $\A$ be a  conformal net. A  conformal net $\B$ on a Hilbert
space $\H$ is an extension of $\A$ or $\A$ is a subnet of $\B$ if
there is a DHR representation $\pi$ of $\A$ on $\H$ such that
$\pi(\A)\subset \B$ is a conformal subnet. The extension is
irreducible if $\pi(\A(I))'\cap \B(I) = {\mathbb C} $ for some (and
hence all) interval $I$, and is of finite index if
$\pi(\A(I))\subset \B(I)$ has finite index for some (and hence all)
interval $I$. The index will be called the index of the inclusion
$\pi(\A)\subset \B$ and is denoted by $[\B:\A].$ If $\pi$ as
representation of
 $\A$ decomposes as $[\pi]= \sum_\lambda m_\lambda[\lambda]$ where
$m_\lambda$  are non-negative  integers and $\lambda$ are
irreducible DHR representations of $\A$, we say that $[\pi]=
\sum_\lambda m_\lambda[\lambda]$ is the spectrum of the extension.
For simplicity we will write  $\pi(\A)\subset \B$ simply as
$\A\subset \B$.
\end{definition}

\subsection{Induction}\label{ind}
Let $\B$ be a conformal net and $\A$ a subnet.
We assume that $\A$ is strongly additive and $\A\subset \B$ has
finite index.
Fix an
interval $I_0\in\I$ and  canonical endomorphism (cf. \cite{LR})
$\gamma$ associated with $\A(I_0)\subset\B(I_0)$.
Given a DHR endomorphism $\r$ of $\B$ localized in $I_0$, the
$\a$-induction $\a_{\r}$
of $\r$ is the endomorphism of $\B(I_0)$ given by
\[
\a_{\r}\equiv \gamma^{-1}\cdot\Ad\e(\r,\l)\cdot\r\cdot\gamma\
\]
where $\e$ denotes the right braiding  (cf. \cite{Xu-m}).

Note that $\Hom( \a_\lambda,\a_\mu)=:\{ x\in \B(I_0) |
x  \a_\lambda(y)= \a_\mu(y)x, \forall y\in \B(I_0)\} $ and
$\Hom( \lambda,\mu)= :\{ x\in \A(I_0) |
x  \lambda(y)= \mu(y)x, \forall y\in \A(I_0)\} .$

%---------------------------------------------------------------------------

\subsection{Preliminaries on VOAs}\label{voa}

We first recall some basic notions from \cite{FLM, Zhu, DLM1}. Let $V =\oplus_{n\geq 0}V_n$ be a vertex
operator algebra as defined in \cite{FLM} (see also \cite{B}). $V$ is called of {\em CFT} type
 if  $\hbox{dim}V_0=1.$ A {\em weak V-module} $M$ is a vector space
equipped with a linear map
\begin{equation*}
\begin{split}
Y_M(\cdot, z) : &V \to (\End M)[[z, z^1]]\\
&v  \mapsto Y_M(v, z) =\sum_{n\in\Z}v_nz^{-n-1}~~~~~~(v_n \in \hbox{End}M)
\end{split}
\end{equation*}
which satisfies the following conditions for $u \in V,$ $v \in V,$ $w \in M$ and $ n \in \Z,$
\begin{equation*}
\begin{split}
&u_nw = 0 \hbox{ for } n >> 0;\\
& Y_M(\1, z) = \hbox{id}_M;\\
z_0^{-1}\delta(&\frac{z_1-z_2}{z_0})Y_M(u, z_1)Y_M(v, z_2)-z_0^{-1}\delta(\frac{z_2-z_1}{-z_0})Y_M(v, z_2)Y_M(u, z_1)\\
&= z_2^{-1}\delta(\frac{z_1-z_0}{z_2})Y_M(Y (u, z_0)v, z_2).
\end{split}
\end{equation*}

A {\em (ordinary) $V$-module} is a weak $V$-module $M$ which carries a $\C$-grading induced by the spectrum
of $L(0)$ where $L(0)$ is a component operator of
$$Y_M(\omega, z) =\sum_{n\in\Z}L(n)z^{-n-2}.$$
That is, $M =\oplus _{\l\in \C}M_{\l}$ where $M_\l = \{w\in M|L(0)w = \l w\}$. Moreover one requires that $M_\l$ is
finite dimensional and for fixed $\l$, $M_{n+\l}= 0$ for all small enough integers $n$.
An {\em admissible} $V$-module is a weak $V$-module $M$ which carries a $\Z_+$-grading
$M =\oplus_{n\in\Z_+} M(n)$ that satisfies the following
$$v_mM(n) \subset M(n+\wt v-m-1)$$
for homogeneous $v\in V.$ It is easy to show that any {\em ordinary} module is {\em admissible}.
And for an {\em admissible} $V$-module $M=\bigoplus_{n\in \Z_+}M(n),$ the contragredient module $M'$
is defined in \cite{FHL} as follows:
\begin{equation*}
M'=\bigoplus_{n\in \Z_+}M(n)^{*},
\end{equation*}
where $M(n)^*=\Hom_{\C}(M(n),\C).$ The vertex operator
$Y_{M'}(a,z)$ is defined for $a\in V$ via
\begin{eqnarray*}
\langle Y_{M'}(a,z)f,u\rangle= \langle f,Y_M(e^{zL(1)}(-z^{-2})^{L(0)}a,z^{-1})u\rangle,
\end{eqnarray*}
where $\langle f,w\rangle=f(w)$ is the natural paring $M'\times M\to \C.$
$V$ is called {\em rational} if every admissible $V$-module is completely reducible.  It is proved in \cite{DLM2}
that if $V$ is rational then there are only finitely many irreducible admissible $V$-modules up to
isomorphism and each irreducible admissible $V$-module is ordinary.  Let
$M^0, \cdots,M^p$ be the irreducible modules up to isomorphism with $M^0 = V$. Then there exist $h_i \in\C$
for $i = 0, \cdots, p$ such that
$$M^i = \oplus_{n=0}^\infty M^i_{h_i+n}$$
where $M^i_{h_i}\neq 0$ and $L(0)|_{M^i_{h_i+n}}= h_i + n,$  $\forall n\in \Z_+.$ $h_i$ is called the {\em conformal weight} of
 $M^i.$  We denote $M^i(n)=M^i_{h_i+n}.$ Moreover, $h_i$ and the central charge $c$ are
rational numbers (see \cite{DLM3}).
Let $h_{\min}$ be the minimum of $h_i$'s. The effective central charge $\tilde c$ is defined as $c-24h_{\min}.$ For each
$M^i$ we define the $q$-character of $M^i$ by
$$\ch_qM^i=q^{-c/24}\sum_{n\geq 0}(\dim M^i_{h_i+n})q^{h_i+n}.$$

$V$ is called $C_2$-{\em cofinite} if $\hbox{dim} V/C_2(V ) <\infty$ where $ C_2(V ) = \la u_{-2}v|u, v \in V \ra$ \cite{Zhu}. Rationality
and $C_2$-cofiniteness are two important concepts in the theory of vertex operator algebras as most
good results in the field need both assumptions.

Take a formal power series in $q$ or a complex function $f(z)=q^h\sum_{n \geq 0} a_nq^n.$ We say that the coefficients of
$f(q)$ satisfy the {\em polynomial growth condition} if there exist positive numbers $A$ and $\alpha$ such that
$|a_n|\leq An^{\alpha}.$

If $V$ is rational and $C_2$-cofinite, then $\ch_qM^i$ converges to a holomorphic function on
the upper half plane \cite{Zhu}. Using the modular invariance result from \cite{Zhu} and results on
vector valued modular forms from \cite{KM} we have (see \cite{DM})

\bl{polygrow}
Let $V$ be rational and $C_2$-cofinite. For each $i,$ the coefficients of $\eta(q)^{\tilde c}\ch_qM^i$ satisfy
 the polynomial growth condition where
 $\eta(q)=q^{1/24}\prod_{n\geq 1}(1-q^n).$
\el

\begin{definition}\label{intertwining}
Let $V$ be a \voa and let $(M^{i},Y_i),$ $(M^{j},Y_j),$ $(M^{k},Y_k)$ be three $V$-modules.
An intertwining operator of type
$\left (
\begin{array}{c}
\ \ M^{k} \ \\
M^{i} \ \  M^{j}
\end{array}
 \right )$ is a linear map
 \begin{eqnarray*}
\begin{aligned}
  \Y (\cdot ,z): &~~~M^{i}&&\to ~~(\Hom(M^{j},M^{k}))\{z\}\\
&v \in M^{i} && \mapsto ~~\Y(v,z)\in (\Hom (M^{j}, M^{k}))\{z\}
\end{aligned}
\end{eqnarray*}
 satisfying the following axioms:

 1. For any $u\in M^{i},$ $$ \Y(L(-1)u,z)=\frac{d}{dz}\Y(u,z),$$
 where $L(n)$ is the component operator of  $Y_{i}(w,z)=\sum_{n\in \Z}L(n)z^{-n-2};$

 2. $\forall u\in V,~v\in M^{i},$
 \begin{eqnarray*}
& z_0^{-1}\d(\frac{z_1-z_2}{z_0})Y_k(u,z_1)\Y(v,z_2)-z_0^{-1}\d(\frac{z_2-z_1}{-z_0})\Y(v,z_2)Y_j(u,z_1)\\
& =z_2^{-1}\d(\frac{z_1-z_0}{z_2})\Y(Y_i(u,z_0)v,z_2).
 \end{eqnarray*}
\end{definition}

The intertwining operators of type $\left (
\begin{array}{c}
\ \ M^{k} \ \\
M^{i} \ \  M^{j}
\end{array}
 \right )$ form a vector space denoted by $\mathcal{V}^{k}_{i,j}.$ The dimension of this vector space
 is called the $fusion~rule$ for $M^{i},$ $M^{j}$ and $M^{k},$ and is denoted by  $N^{k}_{i,j}.$
  We will use $\Y_{i,j}^{k}$ to denote an intertwining operator in $\V^{k}_{i,j}.$
 Assume that $M^s=\sum_{n\in\Z_+} M^s_{\lambda_s+n}$ for $s=i,j,k.$
Then for any $\Y \in \V_{i,j}^{k},$ we know from \cite{FHL} that for $u\in M^{i}$ and $v\in M^{j}$
$$\Y(u,z)v\in z^{\Delta(\Y)}M^{k}[[z,z^{-1}]],$$
where $\Delta(\Y)=\l_{k}-\l_{i}-\l_{j}.$

We now turn our discussion to four point  functions (correlation functions). Let $V$ be a rational
and $C_2$-cofinite vertex operator algebra of {\em CFT} type and $V\cong V'.$ By Lemma 4.1
in \cite{H3}, one knows that  for $u_{a_i}\in M^{a_i},$
 $$\la u_{a'_4}, \Y_{a_1,a_5}^{a_4}(u_{a_1}, z_1)\Y_{a_2,a_3}^{a_5}(u_{a_2},z_2)u_{a_3}\ra,
 $$
$$\la u_{a'_4}, \Y_{a_2,a_6}^{a_4}(u_{a_2}, z_2)\Y_{a_1,a_3}^{a_6}(u_{a_1},z_1)u_{a_3}\ra,
$$
are analytic on $|z_1|>|z_2|>0$ and $|z_2|>|z_1|>0$ respectively and can both be analytically
extended to multi-valued analytic  functions on
$$
R=\{(z_1,z_2)\in \C^2|z_1,z_2\neq 0, z_1\neq z_2\}.$$
We can lift the multi-valued analytic functions on $R$ to single-valued analytic functions on the universal covering
$\tilde R$ of
$R$ as in \cite{H4}. We use
$$E\la u_{a'_4}, \Y_{a_1,a_5}^{a_4}(u_{a_1}, z_1)\Y_{a_2,a_3}^{a_5}(u_{a_2},z_2)u_{a_3}\ra
$$
and
$$E\la u_{a'_4}, \Y_{a_2,a_6}^{a_4}(u_{a_2}, z_2)\Y_{a_1,a_3}^{a_6}(u_{a_1},z_1)u_{a_3}\ra
$$
 to denote those analytic functions.

Let $\{\Y_{a,b;i}^{c}| i=1,\cdots, N_{a,b}^{c}\}$ be a basis of $\V_{a,b}^c.$ The linearly independency of
 $$\{E\la u_{a'_4}, \Y_{a_1,a_5;i}^{a_4}(u_{a_1}, z_1)\Y_{a_2,a_3;j}^{a_5}(u_{a_2},z_2)u_{a_3}\ra|i=1,\cdots,N_{a_1,a_5}^{a_4},j=1,\cdots, N_{a_2,a_3}^{a_5},\forall a_5\}$$
 follows from \cite{H4}.

\subsection{Primary fields for affine VOA and KZ equation}\label{affine}

In this section, we briefly review the construction of the affine \voa associated to the integrable highest weight modules for
the affine Kac-Moody Lie algebras, and also give the KZ equation \cite{KZ,KT} of the correlation functions.

\begin{definition}
Let $W$ be a vector space, a weak vertex operator on $W$ is a formal series
$$a(z)=\sum_{n\in \Z}a_nz^{-n-1}\in (\End ~W)[[z,z^{-1}]]$$
such that for every $w\in W,$ $a_n w=0,$ for $n$ sufficiently large.
\end{definition}

Let $a(z)$ and $b(z)$ be two weak vertex operators on $W,$ define
\be{OPE}
a(z)_nb(z)=\Res_{z_1} ((z_1-z)^na(z_1)b(z)-(-z+z_1)^n b(z)a(z_1)).
\ee
This  is also a weak vertex operator on $W$ (see \cite{LL}).

The following important lemma will be useful later.
\bl{Dong} \cite{Li}
If $a(z),$ $b(z)$ and $c(z)$ are pairwise mutually local weak vertex operators,
 then $a(z)_nb(z)$ and $c(z)$ are mutually local.
\el

Let $\frak{g}$ be a finite-dimensional simple Lie algebra with a
nondegenerate symmetric invariant bilinear form and a Cartan subalgebra $\frak{h}.$
Let  $\hat{\frak{g}}=\C[t, t^{-1}]\otimes \frak{g}\oplus \C K$ be the corresponding
affine Lie algebra. For any $X\in \frak{g},$ set $X(n)=X\otimes t^n$ and
$X(z)=\sum_{n\in \Z}X(n)z^{-n-1}.$   Fix a positive integer $k.$ Then any
$\l \in \frak{h}^*$ can be viewed as a linear form on $\C K\oplus \frak{h}\subset \hat{\frak{g}}$
 by sending $K$ to $k.$ Let us denote the corresponding irreducible
 highest weight module for $\hat{\frak{g}}$ associated to a highest weight $\l$  by $L_{\frak{g}}(k,\l).$
 It is proved that $L_{\frak{g}}(k,0)$ is a rational \voa \cite{DL,FZ,Li}
 with all the inequivalent irreducible modules
 $\{L_{\frak{g}}(k,\l)| \la \l,\theta \ra\leq k, \l\in \frak h^*, \hbox{$\l$ is an integral dominant weight}\} ,$
 where $\theta$ is the longest root of $\frak g$ and $(\theta, \theta)=2.$

$L_{\frak{g}}(k,0)$ has a basis $\{X_{i_1}(-n_1)\cdots X_{i_t}(-n_t)\1| X_{i_s}\in \frak g, n_s \in \Z_+, s=1,\cdots, t \}.$
The vertex operator on $L_{\frak{g}}(k,0)$ is defined as
\begin{eqnarray*}
\begin{aligned}
&Y(X(-1)\1,z)=\sum_{n\in \Z}X(n)z^{-n-1};\\
&Y(X_{i_1}(-n_1)\cdots X_{i_t}(-n_t)\1)=X_{i_1}(z)_{-n_1}\cdots X_{i_t}(z)_{-n_t}\1.
\end{aligned}
\end{eqnarray*}

Let $d=\dim\frak{g},$
 and let $\{u^{(1)},\cdots, u^{(d)}\}$ be an orthogonal basis of $\frak g$ with respect to the bilinear
 form on $\frak g.$ Then set $$\omega = \frac{1}{2(k+\check h)}\sum_{i=1}^d u^{(i)}(-1)u^{(i)}(-1)\1,$$
 where $\check{h}$ is the dual Coxeter number of $\frak{g}.$
Define operators $L(n)$ for $n\in \Z$ by:
$$Y(\omega,z)=\sum_{n\in \Z}L(n)z^{-n-2}.$$
The operators $L(n)$ gives representation of the Virasoro algebra on any $L_{\frak{g}}(k,0)$-modules
with central charge $c=\frac{k\cdot \dim \frak{g}}{2(k+\check{h})}.$
Let $\Y(\cdot,z)$ be an intertwining operator of type
$$\left (
\begin{array}{c}
\ \ L_{\frak{g}} (k,\l_3) \ \\
L_{\frak{g}}(k,\l_2) \ \ L_{\frak{g}}(k,\l_1)
\end{array}
 \right ),$$
 then by \cite{FZ}
  $$
  \Y(u,z)=\sum_{n \in \Z}u_nz^{-n-1}z^{\Delta(\Y)},$$
 where $u_n\in \Hom(L_{\frak{g}}(k,\l_1),L_{\frak{g}}(k,\l_3)),$ and
 $u_n L_{\frak{g}}(k,\l_1)(m)\subset L_{\frak{g}}(k,\l_3)(m+\wt u -n-1),$ where $\wt u=i$  means that
 $u\in L_{\frak{g}}(k,\l_2)(i).$
The following commutator formula  for $u\in L_{\frak g}(k,\l_2)(0)$  is a direct result of the Jacobi identity:
 \begin{equation}\label{comm}
 [X(m),\Y(u,z)]=z^m\Y(X(0)u,z).
 \end{equation}

From now on, we restrict our discussion on affine \voas associated to a finite dimensional simple Lie algebra $\frak{g}$ of level $k$.
By abusing of notation, we use  $\l$ to denote the irreducible $L_{\frak{g}}(k,0)$-module $L_{\frak{g}}(k,\l)$ and $\l'$ to
denote the contragredient module of $\l.$

Let $\l_1,\l_2,\l_3,\l_4$ be four irreducible $L_{\frak{g}}(k,0)$-modules, and fix a basis of intertwining
 operators as in $\S$\ref{voa}. It is proved \cite{KZ, KT} that
\begin{eqnarray}
&span\{E\la u_{\l'_4}, \Y_{\l_3,\mu;i}^{\l_4}(u_{\l_3}, z_1)\Y_{\l_2,\l_1;j}^{\mu}(u_{\l_2},z_2)u_{\l_1}\ra|i,j,\mu\}\\
&=span\{E\la u_{\l'_4}, \Y_{\l_2,\gamma;k}^{\l_4}(u_{\l_2}, z_2)\Y_{\l_3,\l_1;l}^{\gamma}(u_{\l_3},z_1)u_{\l_1}\ra|k,l,\gamma\},
\end{eqnarray}
where $u_{\l_i}\in L_{\frak{g}}(k,\l_i)$.
Then there exist
$(B_{\l_4,\l_1}^{\l_3,\l_2})_{\mu,\gamma}^{i,j;k,l}\in \C$ such that
\be{braiding}
\begin{aligned}
E\la u_{\l'_4}, &\Y_{\l_3,\mu}^{\l_4}(u_{\l_3}, z_1)\Y_{\l_2,\l_1}^{\mu}(u_{\l_2},z_2)u_{\l_1}\ra\\
=& \sum_{k,l,\gamma} (B_{\l_4,\l_1}^{\l_3,\l_2})_{\mu,\gamma}^{i,j;k,l}E
\la u_{\l'_4}, \Y_{\l_2,\gamma;k}^{\l_4}(u_{\l_2}, z_2)\Y_{\l_3,\l_1;l}^{\gamma}(u_{\l_3},z_1)u_{\l_1}\ra,
\end{aligned}
\ee
(see \cite{H2,H3}).  $B_{\l_4,\l_1}^{\l_3,\l_2}$ is called the braiding matrix. In the $\frak{sl}(2)$ case,
since the fusion rule is either $0$ or $1,$ the braiding matrix can be simply
denoted by $(B_{\l_4,\l_1}^{\l_3,\l_2})_{\mu,\gamma},$ since $i,j,k,l \in\{0,1\}.$

Now let us turn our discussion to KZ equations for $L_{\frak{g}}(k,0)$.
 For any
intertwining operators $\Y_2\in \V_{\l_3,\l}^{\l_4}$ and $\Y_1\in \V_{\l_2,\l_1}^{\l},$ and
$u_i\in L_{\l_i}=L_{\frak{g}}(k,\l_i)(0),$ $i=1,2,3,$ $u_4'\in L_{\l'_4}=L_{\frak{g}}(k,\l_4)(0)^*,$ the function
$\Psi(\Y_2,\Y_1,z_1,z_2)$ is defined by:
$$\Psi(\Y_2,\Y_1,z_1,z_2)(u_4\otimes u_3\otimes u_2\otimes u_1)=\la u_4',\Y_2(u_3,z_1)\Y_1(u_2,z_2)u_1 \ra .$$

\bl{} In the region $|z_1|>|z_2|>0,$ this function has a convergent Lauren expansion (\cite{KT}):
$$\Psi(\Y_2,\Y_1,z_1,z_2)(u_4\otimes u_3\otimes u_2\otimes u_1)=
 z_2^{\Delta{(\Y_1)+\Delta(\Y_2)}}\sum_{n\geq 0}\la u_4',(u_3)_{n-1} (u_2)_{-n-1}u_1 (\frac{z_2}{z_1})^{n-\Delta(\Y_2)}\ra.$$
\el
\proof
Direct calculation gives
\begin{equation*}
\begin{aligned}
&\Psi(\Y_2,\Y_1,z_1,z_2)(u_4\otimes u_3\otimes u_2\otimes u_1)\\
=&\la u_4',\Y_2(u_3,z_1)\Y_1(u_2,z_2)u_1 \ra \\
=&\la u_4',\sum_{n,m \in \Z}(u_3)_n (u_2)_m u_1z_1^{-n-1+\Delta(\Y_2)}z_2^{-m-1+\Delta(\Y_1)}  \ra\\
=&\la u_4', \sum_{n\geq 0}(u_3)_{n-1} (u_2)_{-n-1}u_1z_1^{-n+\Delta(\Y_2)}z_2^{n+\Delta(\Y_1)}\ra\\
=& z_1^{\Delta(\Y_2)}z_2^{\Delta(\Y_1)}\sum_{n\geq 0}\la u_4',(u_3)_{n-1} (u_2)_{-n-1}u_1(\frac{z_2}{z_1})^{n}  \ra \\
=&  z_2^{\Delta{(\Y_1)+\Delta(\Y_2)}}\sum_{n\geq 0}\la u_4',(u_3)_{n-1} (u_2)_{-n-1}u_1 (\frac{z_2}{z_1})^{n-\Delta(\Y_2)}\ra.
\end{aligned}
\end{equation*}
\endproof

Introduce a variable $\xi=\frac{z_2}{z_1},$ then the function
$z_2^{-\Delta{(\Y_1)-\Delta(\Y_2)}}\Psi(\Y_2,\Y_1,z_1,\xi z_1)$
is independent of $z_1.$ We abbreviate it to $\Psi(\Y_2,\Y_1,\xi).$
%  \be{red4pt}\psi (\Y_2,\Y_1,\xi)(u_4',u_3,u_2,u_1)=\xi ^{-\Delta(\Y_2)}\sum_{n\geq 0}\la u_4',(u_3)_{n-1} (u_2)_{-n-1}u_1\ra \xi^{n}.\ee

In the case that $u_1\in L_{\l_1}=L_{\frak{g}}(k,\l_1)(0)$ is the highest weight vector of $\frak g$
and $u'_4\in L_{\frak{g}}(k,\l_4)(0)^*$ the lowest weight vector of $\frak g$,
$\Psi(\Y_2,\Y_1,\xi)(u_4',u_3,u_2,u_1)$ verifies the reduced KZ equation \cite{KT}:
 \be{rKZ}
 (k+h^\vee)\frac{d}{d\xi}\Psi =\frac{\Omega_{1,2}-(k+h^{\vee})
 (\Delta(\Y_1)+\Delta(\Y_2))}{\xi}\Psi+\frac{\Omega_{2,3}}{\xi-1}\Psi,
 \ee
 where $\Omega$ is the Casimir element $$\Omega=\sum_{i=1}^d u^{(i)}\otimes u^{(i)},$$
  and
  $$\Omega_{1,2}\Psi (\Y_2,\Y_1,\xi)(u_4',u_3,u_2,u_1)
  =\sum_{i=1}^d \Psi (\Y_2,\Y_1,\xi)(u_4',u_3,a^{(i)}u_2,a^{(i)}u_1),
  $$
  $$\Omega_{2,3}\Psi (\Y_2,\Y_1,\xi)(u_4',u_3,u_2,u_1)
  =\sum_{i=1}^d \Psi (\Y_2,\Y_1,\xi)(u_4',a^{(i)}u_3,a^{(i)}u_2,u_1).$$
  The solutions of the reduced
  KZ equations can only have poles of  finite order.

 \subsection{Crossing symmetry of four point functions}
Crossing symmetry is an important property of quantum groups. Much study was devoted to the connection between conformal field
theory and representations of the braid group (\cite{KT,Ver}, et).
 It is expected that the crossing symmetry of the correlation functions in conformal
 field theory  comes from the crossing symmetry of the quantum group relation.
 The cases of minimal series were carefully treated in \cite{FFK}.
 For the WZW $SU(2)$ model, the elements of braiding matrices  of the correlation functions are essentially the quantum
 $6j$ symbols \cite{Hou}. Following from \cite{AGS,DF1, DF2,Hou}, the crossing symmetry of the braiding matrix
 of the correlation functions in the WZW $SU(2)$ model is derived, and is related to
 the symmetry of the quantum $6j$ symbols.

% In $\S$ \ref{basics}, we discussed the braiding matrices, and we know that the matrices depend on the choices of the basis of intertwining operators. The following Lemma is proved as (3.10) in \cite{}, where the integer level case is treated carefully in $\S$4 (cf. \cite{Felder}) for similar treatment of minimal model case.

 %\begin{lem}\label{crossing symmetry}
%Let $\frak{g}=sl_2(\C)$ and $k\geq 0$ be an integer.  The irreducible modules for $L_{\frak{g}}(k,0)$ are $\{L_{\frak{g}}(k,i)|i=0,\cdots ,k\}.$ The fusion rules are $$i\times j= |i-j|, |i-j|+2, \cdots , \min \{i+j,i+j-k\}.$$
 %\end{lem}

 Now fix a basis of intertwining operators
$\Y_{i,j}^k\in \V_{i,j}^k$ in the $L_{\frak{sl}(2)}(k,0)$ case (here $N_{i,j}^k$ is either $0$ or $1$).
 We use the notation from $\S$\ref{voa} and express the braiding matrix in the following way:
$$
E\la u'_k, \Y_{j,a}^k(u_j,w)\Y_{i,l}^a(u_i,z)u_l  \ra=
\sum_{b}(B^{j,i}_{k,l})_{a,b}E\la u'_k, \Y_{i,b}^k(u_i,z)\Y_{j,l}^b(u_j,w)u_l  \ra.$$

The crossing symmetry can be explained as the following lemma (cf. Page 657 of \cite{FFK}):
\bl{cross sym} The braiding matrix under the basis of intertwining operators
chosen above satisfies
$$(B^{j,i}_{k,l})_{a,b}T_{j,a}^kT_{i,l}^a= (B^{i,j}_{k,l})_{b,a}T_{i,b}^kT_{j,l}^b,$$
where $T_{m,n}^r\in \C$ is a constant uniquely determined by the vertex of type
$\left (
\begin{array}{c}
\ \ r \ \\
m \ \  n
\end{array}
 \right ).$
\el

\subsection{Level-Rank Duality}\label{lrd}
Level-rank duality has been explained by different methods in \cite{GW, SA, NT}.
We will be interested in the following conformal inclusion:
\be{levelrank}
L_{\frak{sl}(m)}(n,0)\otimes L_{\frak{sl}(n)}(m,0)\subset L_{\frak{sl}(mn)}(1,0).
\ee
In the classification of conformal inclusions in \cite{GNO},  the above conformal inclusion
corresponds to AIII.

The decomposition of $L_{\frak{sl}(mn)}(1,0)$ under $L_{\frak{sl}(m)}(n,0)\otimes L_{\frak{sl}(n)}(m,0)$
is known (see \cite{ABI, X-ab}).  To describe such a decomposition, let us prepare some notations. The
level $n$ (resp. $m$) dominant integral weight of $\hat{\frak{sl}}(m)$ (resp. $\hat{\frak{sl}}(n))$ will be denoted by
$\l$ (resp. $\dot \l$). $P_{++}^n$ (resp. $\dot P_{++}^m$)  denote the set of highest weights of level
$n$ of $\hat{\frak{sl}}(m)$ (resp. level $m$ of $\hat{\frak{sl}}(n)$).
The fundamental weight of $\frak{sl}(m)$  (resp. $\frak{sl}(n)$) will be denoted by $\Lambda_i$ (resp. $\dot \Lambda_j$).
We will use $\Lambda_0$ (resp. $\dot \Lambda_0$) or $0$ (resp. $\dot 0$) to denote the trivial representation of $\frak{sl}(m)$ (resp. $\frak{sl}(n)$). Then any $\l$ can be expressed as  $\l=\sum_{i=0}^{m-1}\l_i\L_i,$ and $\sum_{i=0}^{m-1}\l_i=n.$
Instead of $\l=(\l_0,\cdots ,\l_{m-1}),$  it will be more convenient to use
$$
\l+\rho=\sum_{i=0}^{m-1}\l'_i\Lambda_i$$
with $\l'_i=\l_i+1.$ Then
$\sum_{i=0}^{m-1}\l'_i=m+n.$

Due to the cyclic symmetry of the extended Dynkin diagram of $\frak{sl}(m),$ the group $\Z_m$ acts on $P_{++}^n$ by
$$
\Lambda_i\to \Lambda_{(i+\mu)modm},~ \mu\in\Z_m.$$
Let $\Omega_{m,n}=P_{++}^n/\Z_m.$ Then there is a natural bijection between
$\Omega_{m,n}$ and $\Omega_{n,m}$ (see $\S$2 of \cite{ABI}).

We parameterize the bijection by a map
$$
\beta: P_{++}^n\to\dot P_{++}^m$$
as follows. Set
$$
r_j=\sum_{i=j}^m \l'_i,~~ 1\leq j\leq m,$$
where $\l'_m\equiv \l'_0.$ The sequence $(r_1,\cdots, r_m)$ is decreasing, $m+n=r_1>r_2>\cdots >r_m\geq 1.$
Take the complementary sequence $(\bar r_1,\bar r_2,\cdots, \bar r_n)$ in $\{1,2, \cdots, m+n\}$ with
$\bar r_1> \bar r_2 >\cdots > \bar r_n.$
Put $$s_j=m+n+\bar r_n-\bar r_{n-j+1}, ~ 1\leq j\leq n.$$
Then $m+n=s_1>s_2>\cdots>s_n\geq 1.$
The map $\beta$ is defined by
$$
(r_1,\cdots, r_m)\to (s_1,\cdots ,s_n).$$

The following lemmas summarizes what we will use.

\bl{duality}\cite{X-ab} Let $Q$ be the root lattice of $\frak{sl}(m),$
$\Lambda_i,$ $0\leq i\leq m-1,$ its fundamental weights,
and $Q_i=(Q+\Lambda_i)\cap P_{++}^n.$ Let $\tilde \Lambda\in\Z_{mn}$
denote a level 1 highest weight of $\frak{sl}(mn)$
and $\l\in Q_{\tilde{\Lambda} modm}.$ Then there exists a unique
$\dot \l\in \dot P_{++}^m$ with $\dot \l=\mu \beta(\l)$ for some unique $\mu\in \Z_n$
such that
$L_{\frak{sl}(m)}(n,\l)\otimes L_{\frak{sl}(n)}(m,\dot\l)$ appears once
and only once in $L_{\frak{sl}(mn)}(1,\tilde\Lambda).$
The map $$\l \to \dot \l=\mu\beta(\l)$$ is one-to-one. Moreover,
$L_{\frak{sl}(mn)}(1,\tilde\Lambda)$ is a direct sum of all
$L_{\frak{sl}(m)}(n,\l)\otimes L_{\frak{sl}(n)}(m,\dot\l):$
$$
L_{\frak{sl}(mn)}(1,\tilde\Lambda)=\bigoplus_{\l\in Q_{\tilde\Lambda modm}}
L_{\frak{sl}(m)}(n,\l)\otimes L_{\frak{sl}(n)}(m,\dot\l).$$
\el

\bl{hom}
Take $\tilde\Lambda=0$ in the above lemma,  $\{\l|\l\in Q_0\}$ (resp. $\{\dot\l\}$) is closed under fusion.
Moreover, the map  $\l \to \dot\l$ gives an isomorphism between the two fusion subalgebras.
\el

\begin{rmk}\label{branchingrule}
In the case $m=2,$ $n=10$ and $\tilde\Lambda=0,$ $Q_0=\{0,2,4,6,8,10\}.$ When we take $m=2,$ $n=28,$ $Q_0=\{0,2,4,\cdots, 28\}.$

\end{rmk}

%---------------------------------------------------------------------------------
%---------------------------------------------------------------------------------

We write  the conformal net and subnet which correspond to $L_{\frak{sl}(m)}(n,0)\otimes L_{\frak{sl}(n)}(m,0)\subset L_{\frak{sl}(mn)}(1,0)$  as $\A\subset \B.$ For simplicity we assume that the spectrum of $\A\subset \B$ is $\sum_{\lambda} \lambda\otimes (1,\lambda)$ where $\lambda$, $(1,\lambda)$ label the irreducible representations of $L_{\frak{sl}(m)}(n,0)$ and $L_{\frak{sl}(n)}(m,0)$
respectively.

\begin{lem}\label{above}
For $U_{\bar{\l}}\in \Hom(\a_{\bar{\l}},\a_{(1,\l)})$ be a unitary as in (3) of Prop. 3.7 in \cite{Xu-m},  $T_\l\in \Hom
(1,\a_{\l,(1,\l)}),$ one has $U_{\bar{\l}}=n_\l T_\l,$
for some $n_\l \in \A(I).$
\end{lem}
\proof  Let $r_\l \neq 0,$  $r_\l \in \Hom(\l \bar{\l},1).$ Then $(1,\l)
(r_\l)T_\l \in \Hom(\a_{\bar{\l}},\a_{(1,\l)}),$
 and $(1,\l)(r_\l)T_\l \neq 0.$ It follows that $U_{\bar{\l}}=n_\l T_\l
 ,$ for some $n_\l \in \A(I).$
\endproof

\begin{prop}\label{nond}
If $\l_3 \prec \l_1 \l_2,$ then $E(T_{\l_3}T_{\l_1}^*T_{\l_2}^*)\neq 0,$ where $T_{\lambda_i}\in \Hom(1,\a_{\lambda_i}), i=1,2,3.$
\end{prop}

\proof Choose $T\neq 0,$
$T\in \Hom(\overline{\l_3},\overline{\l_1}\overline{\l_2}).$ Then $U_{\l_2}U_{\l_1}
T U_{\l_3}^*\in \Hom(\a_{(1,\l_3)},\a_{(1,\l_1)(\a,\l_2)})= \Hom ((1,\l_3),
(1,\l_1)(1,\l_2))\subset \A(I).$ So $E(U_{\l_2}U_{\l_1}T U_{\l_3}^*)=
U_{\l_2}U_{\l_1}T U_{\l_3}^*\neq 0.$

By the Lemma \ref{above},
\begin{eqnarray*}
\
E(U_{\l_2}U_{\l_1}T U_{\l_3}^*)&&=E(n_{\l_2}T_{\l_2}n_{\l_1}T_{\l_1}
T T_{\l_3}^*n_{\l_3}^*)\\
&&=E(n_{\l_2}\l_2(n_{\l_1})T_{\l_2}T_{\l_1}T_{\l_3}^*\l_3(T)
n_{\l_3}^*)\\
&&=n_{\l_2}\l_2(n_{\l_1})E(T_{\l_2}T_{\l_1}T_{\l_3}^*)\l_{3}(T)n_{\l_3}^*)\\
&&\neq 0.
\end{eqnarray*}
It follows that $E(T_{\l_2}T_{\l_1}T_{\l_3}^*) \neq 0.$
Taking the adjoint we have proved the Proposition.
\endproof

\begin{cor}\label{nondeg}
For  inclusions $SU(n)_m\times SU(m)_n\subset SU(mn)_1,$
suppose the vertex operator
$$J^{(\l_2,\dot{\l_2})}(z)=\sum_{\l_1,\l_3}D_{\l_2,\l_1}^{\l_3}
\Y_{\l_2,\l_1}^{\l_3}(\cdot ,z)\otimes \Y_{\dot{\l}_2,\dot{\l}_1}^{\dot{\l}_3}
(\cdot ,z).$$
If $\l_3\prec \l_2 \l_1 ,$ then $D_{\l_2,\l_1}^{\l_3}\neq 0. $
\end{cor}
\proof  If $D_{\l_2,\l_1}^{\l_3}=0,$ then $H_{\l_3}\otimes H_{\dot{\l}_3}
\bot J^{(\l_2,\dot{\l_2})}J^{({\l}_1,\dot{\l}_1)}H_0\otimes H_{\dot{0}}.$
By Proposition \ref{nond} and (2) of Lemma 3.3 in \cite{Xeg}, we have
$H_{\l_3}\otimes H_{\dot{\l}_3}\subset J^{(\l_2,\dot{\l}_2)}J^{(\l_1,\dot{\l}_1)}
H_0\otimes H_{\dot{0}},$ where $J^{(\l_2,\dot{\l}_2)}=V(\l_2,\dot{\l}_2)$ in (2) of Lemma 3.3 in \cite{Xeg},  a contradiction.
\endproof

\begin{rmk}
A proof of Corollary \ref{nondeg} using \voa  language has not been found in this paper. We will do direct calculation in the case we need for constructing the main examples (see Remark \ref{voa nondeg} and $\S$\ref{appendix}).
\end{rmk}

\subsection{Lattice Vertex Operator Algebras}\label{lattice}
%In this section we recall the construction of the \voa $V_L$ associated with a positive definite even lattive $L$ and its irreducible modules.
Let $L$ be a rank $d$ even lattice with a positive definite symmetric $\Z$-bilinear form
$(\cdot ,\cdot).$ We set $\frak{h}=\C \otimes_\Z L$ and extend $(\cdot , \cdot)$ to a $\C$-bilinear
form on $\frak{h}.$ Let $\hat{\frak{h}}=\C[t,t^{-1}]\otimes \frak{h}\oplus \C C$ be the affinization
of commutative Lie algebra ${\frak{h}}.$ For any $\l \in \frak{h},$ we can define a one
dimensional $\hat{\frak{h}}^+$-module $\C e^\l $ by the actions
$\rho(h\otimes t^m)e^\l=(\l,h)\d_{m,0}e^\l$ and $\rho (C) e^\l=e^\l $ for $h\in \frak{h}$ and
$m\geq 0.$ Now we denote by
$$M(1,\l)=U(\hat{\frak{h}})\otimes_{U(\hat{\frak{h}}^+)}\C e^\l\cong S(t^{-1}\C[t^{-1}])$$
the $\hat{\frak{h}}$-module induced from $\hat{\frak{h}}^+$-module. Set $M(1)=M(1,0).$
Then there exists a linear map $Y: M(1)\to (\End M(1,\l))[[z,z^{-1}]]$ such that $(M(1), Y, \1, \o)$
is a simple \voa  and $(M(1,\l),Y)$ becomes an irreducible $M(1)$-module for any
$\l \in \frak{h}$ (see \cite{FLM}). The vacuum vector and the Virasoro element are given by $\1 = e^0$ and
$\o =\frac{1}{2}\sum_{i=1}^d a_i(-1)^2\otimes e^0,$ respectively, where $\{a_i\}$ is an
orthonormal basis of $\frak{h}.$

Let $L$ be any positive definite even lattice and let $\hat{L}$ be the canonical central
extension of $L$ by the cyclic group $\la \kappa \ra$ of order 2:
$$1\to \la \kappa \ra \to \hat{L}\bar{\to}L\to 0$$
with the commutator map
$c(\a ,\b)=\kappa ^{(\a , \b)}$ for $\a ,\b \in L.$ Let $e:L\to \hat{L}$
be a section such that
$e_0=1$ and $\varepsilon :L\times L\to \la \kappa \ra$ be the
corresponding 2-cocycle. We can assume that $\varepsilon$ is bimultiplicative. Then
$\varepsilon(\a ,\b)\varepsilon(\b,\a)=\kappa^{(\a ,\b)},$
$$\varepsilon(\a ,\b)\varepsilon(\a +\b ,\g)=\varepsilon(\b ,\g)\varepsilon(\a ,\b+\g)$$
 and $e_\a e_\b=\varepsilon(\a ,\b)e_{\a +\b}$ for $\a ,\b ,\g \in L.$

Let $L^{\circ} =\{\l\in \frak{h}|(\a ,\l)\in \Z\}$ be the dual lattice of $L.$ Then there is a
$\hat{L}$-module structure on $\C[L^{\circ}]=\bigoplus_{\l \in L^\circ}\C e^\l$ such that
$\kappa$ acts as $-1$ (see \cite{DL}). Let $L^\circ=\bigcup_{i\in L^\circ /L}(L+\l_i)$ be the coset
decomposition such that $\l_0=0.$ Set $\C[L+\l_i]=\bigoplus_{\a \in L}\C e^{\a+\l_i}.$
Then $\C[L^\circ]=\bigoplus_{i\in L^\circ /L}\C[L+\l_i]$ and each  $\C[L+\l_i]$ is an
 $\hat{L}$-submodule of $\C[L^\circ].$ The action of $\hat{L}$ on $\C[L+\l_i]$ is as follows:
$$e_\a e^{\b+\l_i}=\varepsilon(\a ,\b)e^{\a +\b +\l_i}$$
for $\a ,\b\in L.$

We can identify $e^\a$ with $e_\a$ for $\a \in L.$ For any $\l \in L^\circ,$ set
$\C[L+\l]=\bigoplus_{\a \in L}\C e^{\a+\l}$ and define $V_{L+\l}=M(1)\otimes \C[L+\l].$
Then for any $V_{L+\l},$
there exists a linear map
$$Y: V_L \to (\End V_{L_\l})[[z,z^{-1}]]$$
such that $(V_L,Y,\a,\o)$
becomes a simple \voa and $(V_{L+\l},Y)$ is an irreducible $V_L$-module \cite{B, FLM}.
And $V_{L+\l_i}$ for $\l_i \in L^\circ /L$ give all inequivalent irreducible $V_L$-modules (see \cite{D}).
The vertex operator $Y(h(-1)\1,z)$ and $Y(e^\a ,z)$ associated to $h(-1)\1$ and $e^\a$ are defined as
$$Y(h(-1)\1,z)=h(z)=\sum_{n\in \Z}h(-n)z^{-n-1},
$$
$$Y(e^\a,z)=\exp(\sum_{n=1}^\infty \frac{\a(-n)}{n}z^{n})\exp(-\sum_{n=1}^\infty \frac{\a(n)}{n}z^{-n})e_\a z^{\a},
$$
where $h(-n)$ is the action of $h\otimes t^n$ on $V_{L+\l},$ $e_\a$ is the left action of
$\hat{L}$ on $\C[L^\circ],$ and $z^\a$ is the operator on $\C[L^\circ]$ defined by
 $z^\a e^\l=z^{(\a, \l)}e^\l .$
The vertex operator associated to the vector $v=\b_1(-n_1)\cdots \b_r(-n_r)e^\a$ for
$\b_i \in \frak{h},$ $n_i\geq 1,$ and $\a \in L$ is defined as
$$Y(v,z)=:\partial^{(n_1-1)}\b_1(z) \cdots \partial^{(n_r-1)}\b_r(z)Y(e^\a,z):,$$
where $\partial=(1/n!)(d/dz)$ and  $:,:$ is the normal ordered products.

In the case we choose $L$ be the root lattice of simple Lie algebras $\frak{g}$ of ADE type,
one knows $V_L\cong L_{\frak{g}}(1,0)$ as vertex operator algebras. %When we consider our key example, we
%need to use the isomorphism $V_{L_{\frak{sl}(20)}\cong L_{sl(20)}(1,0).$

%----------------------------------------------------------------------------------
%----------------------------------------------------------------------------------

\section{The mirror extension of $L_{\frak{sl}(10)}(2,0)$}\label{key}

This section is devoted to the construction of \voa
$$
V=L_{\frak{sl}(10)}(2,0)^{e}=L_{\frak{sl}(10)}(2,0)\oplus L_{\frak{sl}(10)}(2,\Lambda_3+\Lambda_7)$$
based on the conformal inclusions
\be{b2}
SU(2)_{10}\subset Spin(5)_1
 \ \ ~~(L_{\frak{sl}(2)}(10,0)\subset L_{B_2}(1,0))
 \ee
 and
 \be{sl20}
 SU(2)_{10}\times SU(10)_2\subset SU(20)_1\ \ ~~(L_{\frak{sl}(2)}(10,0)\otimes L_{\frak{sl}(10)}(2,0)\subset L_{\frak{sl}(20)}(1,0)).
 \ee

\subsection{The conformal inclusions}

The conformal inclusion (\ref{b2}) is well studied in conformal nets theory.
Due to \cite{CIZ} and  a recent result of \cite{DLN}, the corresponding conformal inclusion
of \voas and the branching rules are also established in \voa theory.
The decomposition of $L_{B_2}(1,0)$ as an $L_{\frak{sl}(2)}(10,0)$-module is as follow:
$$
L_{B_2}(1,0)=L_{\frak{sl}(2)}(10,0)\oplus L_{\frak{sl}(2)}(10,6).
$$
For convenience, we denote the above decomposition as
\be{0+6}L_{B_2}(1,0)=0+6.
\ee

The conformal inclusion (\ref{sl20}) comes from the level-rank duality ($\S$\ref{lrd})
and the branching rules are given in Lemma \ref{duality} and Remark \ref{branchingrule}:
\begin{equation}\label{eqsl20}
\begin{aligned}
L_{\frak{sl}(20)}(1,0))&=  L_{\frak{sl}(2)}(10,0) \otimes L_{\frak{sl}(10)}(2,0)\\
&\oplus L_{\frak{sl}(2)}(10,2) \otimes L_{\frak{sl}(10)}(2,\Lambda_1+\Lambda_9)\\
&\oplus L_{\frak{sl}(2)}(10,4) \otimes L_{\frak{sl}(10)}(2,\Lambda_2+\Lambda_8)\\
&\oplus L_{\frak{sl}(2)}(10,6) \otimes L_{\frak{sl}(10)}(2,\Lambda_3+\Lambda_7)\\
&\oplus L_{\frak{sl}(2)}(10,8) \otimes L_{\frak{sl}(10)}(2,\Lambda_4+\Lambda_6)\\
&\oplus L_{\frak{sl}(2)}(10,10) \otimes L_{\frak{sl}(10)}(2,2\Lambda_5)
\end{aligned}
\end{equation}
We will use the notation as in Lemma \ref{duality} for the above decomposition:
$$
L_{\frak{sl}(20)}(1,0))=\sum_{\l \hbox{ even},\l=0}^{10} \l \times \dot \l .$$
The decompositions (\ref{0+6}) and (\ref{eqsl20}) and the Mirror Extension Conjecture
allow us to make the following assertion, which is the first main theorem of this paper:
\bt{main}
There is a \voa structure on
$$
V=L_{\frak{sl}(10)}(2,0)^e=\dot 0+\dot6=L_{\frak{sl}(10)}(2,0)\oplus L_{\frak{sl}(10)}(2,\Lambda_3+\Lambda_7).$$
\et

 \subsection{The construction of the VOA extension}

 In this section, we focus on defining the vertex operator
 $$
 \dot Y(\cdot,z) : V\to \End V[[z,z^{-1}]]$$
 that gives a \voa structure on $L_{\frak{sl}(10)}(2,0)^e.$

First we use $Y(\cdot,z)$ and $\tilde Y(\cdot,z)$ to denote the vertex
operators of $L_{B_2}(1,0)$ and $L_{\frak{sl}(20)}(1,0)$ respectively.
We fix a basis of the intertwining operators $\Y_{a,b}^c$
(resp. $\Y_{\dot a,\dot b}^{\dot c}$)  among modules $\{\l\in Q_0\}$ (resp. $\{\dot \l\}$)
of $L_{\frak{sl}(2)}(10,0)$  (resp. $L_{\frak{sl}(10)}(2,0)$) such that
\be{VO SU(20)}
\tilde Y(u_1\otimes u_2,z)=\sum_{\l_1, \l_2\in Q_0}D_{\l,\l_1}^{\l_2}
\mathcal{\Y}^{\l_2}_{\l,\l_1}(u_1,z)\otimes \mathcal{Y}^{\dot{\l_2}}_{\dot{\l},\dot{\l_1}}(u_2,z),
\ee
for $u_1\in L_{\frak{sl}(2)}(10,\l),$ $u_2\in L_{\frak{sl}(10)}(2,\dot\l),$
$\Y^{\l_2}_{\l,\l_1}\in\V^{\l_2}_{\l,\l_1}.$
One can choose $D_{\l,\l_1}^{\l_2}=\delta_{1,N^{\l_2}_{\l,\l_1}}$ because of Corollary \ref{nondeg} by
suitably choosing the basis of intertwining operators $\Y_{a,b}^c$ and $\Y_{\dot a,\dot b}^{\dot c}$.
  But we want a pure vertex operator algebra proof that $D_{\l,\l_1}^{\l_2}\ne 0$
if $N_{\l,\l_1}^{\l_2}\ne 0.$ It turns out that  we only need to show that $D_{6,6}^{\mu}\ne 0,$
$D_{6,\mu}^6\ne 0$ for $\mu=0,2,4,6,8$ for the purpose of this paper.  This result is given in $\S$\ref{appendix}.
So in the discussion below we always assume that $D_{6,6}^{\mu}=D_{6,\mu}^6=D_{\l,\l_1}^{\l_2}=1$ with $N_{\l,\l_1}^{\l_2}\ne 0$ for
$\mu=0,2,4,6,8$ and one of $\l,\l_1,\l_2$ being $0.$

Under the same basis of intertwining operators,
for $u\in L_{\frak{sl}(2)}(10,0)\subset L_{\frak{sl}(2)}(10,0)^e$ the vertex operator is obviously of the form
 $$
Y(u,z)= J^0(u,z)=Y_0(u,z)+Y_6(u,z).$$
 For $u\in L_{\frak{sl}(2)}(10,6),$ we denote $Y(u,z)$ by
$$J^6(u,z)=Y(u,z)=\sum_{\l_1,\l_2\in\{0,6\}} c_{6,\l_1}^{\l_2} \mathcal{Y}_{6,\l_1}^{\l_2}(u,z),$$
 where $c_{6,\l_1}^{\l_2} \in\C.$
 Similarly,  we write
 $$
 \dot J^{\dot0}(u,z)=\dot Y(u,z)=Y_{\dot0}(u,z)+Y_{\dot6}(u,z)$$
for $u\in  L_{\frak{sl}(10)}(2,\dot 0)$ and
$$\dot J^{\dot6}(u,z)=\dot Y(u,z)=\sum_{\l_1,\l_2\in\{0,6\}} c_{\dot6,\dot\l_1}^{\dot\l_2}
\mathcal{Y}_{\dot6,\dot\l_1}^{\dot\l_2}(u,z)
 $$
 for $u\in L_{\frak{sl}(10)}(2,\Lambda_3+\Lambda_7),$
 where $c_{\dot6,\dot\l_1}^{\dot\l_2}\in\C$ are the coefficients needed to be
  determined such that $\dot Y$ satisfying locality.

%where  $\V^{\l_1}_{\l,\l_3}$ is the space of intertwining operators of type $\left (
%\begin{array}{c}
%\ \ \l_3 \ \\
%\l \ \  \l_1
%\end{array}
% \right ),$

Notice that $L_{B_2}(1,0)$ is self dual.
Then for any $u_i \in L_{B_2}(1,0),$ $i=1,2,3,4,$
$$E\langle u_4, Y(u_3,w)Y(u_2,z)u_1 \rangle$$
is a rational symmetric function since $L_{B_2}(1,0)$ is a vertex operator algebra. If we choose
$$u_i\in L_{\frak{sl}(2)}(10,\l_i)\subset L_{B_2}(1,0),  ~~ i=1,2,3,4, $$  where $\l_i\in \{0,6\},$ then we have
\begin{eqnarray}\label{1}
\begin{aligned}
&E\langle u_4, Y(u_3,w)Y(u_2,z)u_1 \rangle\\
=&E\langle u_4, J^{\l_3}(u_3,w)J^{\l_2}(u_2,z)u_1 \rangle\\
=& \sum_{\mu \in \{0,6\}} c_{\l_3,\mu}^{\l_4}c_{\l_2,\l_1}^{\mu}
E\langle u_4, \mathcal{Y}_{\l_3,\mu}^{\l_4}(u_3,w) \mathcal{Y}_{\l_2,\l_1}^{\mu}(u_2,z)u_1 \rangle\\
=&  \sum_{\mu ,\gamma  \in \{0,6\}}  c_{\l_3,\mu}^{\l_4}c_{\l_2,\l_1}^{\mu}(B^{\l_3,\l_2}_{\l_4,\l_1})_{\mu,\gamma}E\langle u_4,
\mathcal{Y}_{\l_2,\gamma}^{\l_4}(u_2,z) \mathcal{Y}_{\l_3,\l_1}^{\gamma}(u_3,w)u_1 \rangle\\
= &\sum_{\gamma \in \{0,6\}}  c_{\l_2,\gamma}^{\l_4}c_{\l_3,\l_1}^{\gamma}E\langle u_4,
\mathcal{Y}_{\l_2,\gamma}^{\l_4}(u_2,z) \mathcal{Y}_{\l_3,\l_1}^{\gamma}(u_3,w)u_1 \rangle\\
=&E\langle u_4, J^{\l_2}(u_2,z)J^{\l_3}(u_3,w)u_1 \rangle\\
=&E\langle u_4, Y(u_2,z)Y(u_3,w)u_1 \rangle
\end{aligned}
\end{eqnarray}
by using the braiding isomorphism and locality of correlation functions of $L_{B_2}(1,0).$

Due to the linearly independent property of the correlation functions (see $\S$\ref{voa}),  we have
\begin{eqnarray}\label{2}
\Sigma_{\mu  \in \{0,6\}}  c_{\l_3,\mu}^{\l_4}c_{\l_2,\l_1}^{\mu}(B^{\l_3,\l_2}_{\l_4,\l_1})_{\mu,\gamma}=  c_{\l_2,\gamma}^{\l_4}c_{\l_3,\l_1}^{\gamma}
\end{eqnarray}

\begin{lem}\label{rational}
 Equation (\ref{2}) is a necessary and sufficient condition for
$$E\langle u_4, J^{\l_3}(u_3,w)J^{\l_2}(u_2,z)u_1 \rangle$$
 to be a symmetric rational function of $w, z,$  for primary vectors $u_i\in L_{\frak{sl}(2)}(10,\l_i)(0),$
 $i=2,3,$ $u_1$ the highest weight vector of $\frak{sl}(2)$ in $L_{\frak{sl}(2)}(10,\l_1)(0)$ and $u_4$ the lowest weight vector of $\frak{sl}(2)$ in
 $L_{\frak{sl}(2)}(10,\l_4)(0)^*,$ where $$J^{\l_i}=\sum_{\l_j,\l_k}c_{\l_i,\l_j}^{\l_k}\Y_{\l_i,\l_j}^{\l_k}.$$
 \end{lem}
\pf
According to equation (\ref{1}), the condition is clearly necessary. Now assume equation (\ref{2}) holds.  Then
$E\langle u_4, J^{\l_3}(u_3,w)J^{\l_2}(u_2,z)u_1 \rangle$ is obviously symmetric.

Since
$$
\langle u_4, J^{\l_3}(u_3,w)J^{\l_2}(u_2,z)u_1 \rangle = z^{h_{\l_4}-h_{\l_1}-h_{\l_2}-h_{\l_3}}\Psi (\xi), $$
 where $\xi =\frac{z}{w}.$ The $\Psi (\xi)$ satisfies the reduced KZ equation (\ref{rKZ}). Thus $\Psi (\xi)$
is analytic except at $0, $ $1,$  $\infty,$ and  can only have
 poles of finite order at $0,$ $1,$ $\infty .$ It follows that $\Psi (\xi)$ is a
 rational function of $\xi .$ %The symmetric follows directly from equation (\ref{1}) without the assumption $u_i$ is primary.
 \qed

In order to prove Theorem \ref{main}, we only need
to define the vertex operator on $V$ satisfying locality, which
is equivalent to that the four point functions are rational symmetric
functions. So it suffices to find solutions to the dotted version of equation (\ref{2}).

The following lemma which is essential in our proof of
 Theorem \ref{main} comes from the locality of vertex operators $\tilde Y$
  (see equation (\ref{VO SU(20)})) on $L_{\frak{sl}(20)}(1,0).$
\begin{lem} We have
\be{eq2}
\sum_{\mu} (B^{\l_3,\l_2}_{\l_4,\l_1})_{\mu,\gamma}\cdot
(B^{\dot\l_3,\dot\l_2}_{\dot{\l_4},\dot{\l_1}})_{\dot{\mu},\dot{\gamma_1}}=\delta_{\gamma,\gamma_1}
\ee
for $\l_i\in \{0,6\},$ $i=1,2,3,4.$
 \end{lem}
 \proof  Since every $L_{\frak{sl}(2)}(10,\l_i)\otimes L_{\frak{sl}(10)}(2,\dot {\l_i})$ is a self dual
 $ L_{\frak{sl}(2)}(10,0)\otimes L_{\frak{sl}(10)}(2,\dot 0)$-module, for $u_i\otimes \dot{u}_i \in L_{\frak{sl}(2)}(10,\l_i)\otimes L_{\frak{sl}(10)}(2,\dot {\l_i}),$
 $i=1,2,3,4,$ %and $u_4\otimes \dot u_4\in L_{\frak{sl}(2)}(10,\l_4')\otimes L_{\frak{sl}(10)}(2,\dot \l_4'),$
  the four point function
\begin{eqnarray}
\begin{aligned}
& E \langle u_4\otimes \dot {u}_4, \tilde Y(u_3\otimes \dot{u}_3,w)\tilde Y(u_2\otimes \dot{u}_2,z)u_1\otimes \dot{u}_1 \rangle\\
=& \Sigma_{\mu} E\langle u_4, \mathcal{Y}_{\l_3,\mu}^{\l_4}(u_3,w) \mathcal{Y}_{\l_2,\l_1}^{\mu}(u_2,z)u_1\rangle
\otimes  \langle \dot{u}_4, \mathcal{Y}_{\dot{\l_3},\dot{\mu}}^{\dot{\l_4}}(\dot{u}_3,w)
\mathcal{Y}_{\dot{\l_2},\dot{\l_1}}^{\dot{\mu}}(\dot{u}_2,z)\dot{u}_1\rangle
\end{aligned}
\end{eqnarray}
is a rational symmetric function.

Since switching  $\l_2, \l_3 $, and $\dot{\l_2}, \dot{\l_3}$ gives the same analytic continuation,  we have
\begin{eqnarray}
\begin{aligned}
\Sigma_{\mu,\gamma,\dot{\gamma}_1}E\la u_4\otimes \dot u_4,& (B^{\l_3,\l_2}_{\l_4,\l_1})_{\mu,\gamma}
\cdot (B^{\dot{\l_3},\dot{\l_2}}_{\dot{\l_4},\dot{\l_1}})_{\dot{\mu},\dot{\gamma_1}}\\
&\mathcal{Y}_{\l_2,\gamma}^{\l_4}(u_2,z) \mathcal{Y}_{\l_2,\l_1}^{\gamma}(u_3,w)
\otimes  \mathcal{Y}_{\dot{\l_2},\dot{\gamma_1}}^{\dot{\l_4}}(\dot{u}_2,z)
\mathcal{Y}_{\dot{\l_3},\dot{\l_1}}^{\dot{\gamma_1}}(\dot{u}_3,w)u_1\otimes\dot u_1\ra \\
=\Sigma_{\gamma}E\la u_4\otimes \dot u_4, &\mathcal{Y}_{\l_2,\gamma}^{\l_4}(u_2,z) \mathcal{Y}_{\l_3,\l_1}^{\gamma}(u_3,w)
\otimes  \mathcal{Y}_{\dot{\l_2},\dot{\gamma}}^{\dot{\l_4}}(\dot{u}_2,z)
 \mathcal{Y}_{\dot{\l_3},\dot{\l_1}}^{\dot{\gamma}}(\dot{u}_3,w)u_1\otimes\dot u_1\ra.
\end{aligned}
\end{eqnarray}
Due to the linear independence of the four point functions for $L_{\frak{sl}(2)}(10,0)\otimes L_{\frak{sl}(10)}(2,0)$
in the above equation, one must have
\begin{eqnarray}\label{3}
\Sigma_{\mu} (B^{\l_3,\l_2}_{\l_4,\l_1})_{\mu,\gamma}\cdot
(B^{\dot{\l_3},\dot{\l_2}}_{\dot{\l_4},\dot{\l_1}})_{\dot{\mu},\dot{\gamma_1}}=\delta_{\gamma,\dot{\gamma}_1}
\end{eqnarray}
as desired. \endproof

We are now in a position  to determine  the vertex operator on
$V=L_{\frak{sl}(10)}(2,0)^e.$
The locality of $\dot Y(u,w)\dot Y(v,z)$ is easy to see for either
$u\in L_{\frak{sl}(10)}(2, 0)$ or $v\in L_{\frak{sl}(10)}(2,0).$
Thus we only need to consider locality of the case
$\dot Y(u,w)\dot Y(v,z)$  for both $\dot u,\dot v\in L_{\frak{sl}(10)}(2,\Lambda_3+\Lambda_7).$
We focus on the four point function
$$
E\la \dot u_4,\dot J^{\dot6}(\dot u_3,w)\dot J^{\dot6}(\dot u_2,z)\dot u_1\ra$$
which need to be rational and
symmetric.

We first only focus on choosing $\dot u_1, \dot u_2,\dot u_3, \dot u_4$  to be primary vectors as in Lemma \ref{rational}.
%Need to find $J^{\dot6}(u,z)=\dot c_{\dot 6,\do0}\Y_{\dot6,\dot0}^{\dot6}
%+\dot c_{\dot6,\dot6}\Y_{\dot6,\dot6}^{\dot6}
%+\dot c_{\dot0,\dot6}\Y_{\dot0,\dot6}^{\dot6}, $ satisfying the dotted equation (\ref{eq1}).
In this case, we only need to check the dotted equation (\ref{2}).

\bl{case1} If one of $\dot\l_1, \dot\l_4$ is $\dot0$, the dotted equation (\ref{2}) is automatically satisfied
(dotted (\ref{2}) is trivial in this case).
\el
\proof  Suppose $\l_4=0, \l_1=6.$ Then $\gamma=\mu=6,$ i.e there is only one possible channel.
Equation (\ref{2}) implies
 $$c_{6,6}^0c_{6,6}^6(B^{6,6}_{0,6})_{6,6}=c_{6,6}^0c_{6,6}^6.$$
 Since $c_{6,6}^0\cdot c_{6,6}^6\neq 0$, one gets $(B^{6,6}_{0,6})=1.$

From equation (\ref{eq2}), we have
$$(B^{6,6}_{0,6})_{6,6}(\dot B^{\dot6,\dot6}_{\dot0,\dot6})_{\dot6,\dot6}=1.$$
which implies $\dot B^{\dot6,\dot6}_{\dot0,\dot6}=1,$ so the dotted equation (\ref{2}) is trivial.
\endproof

 It remains to deal with the case $\dot\l_1=\dot\l_2=\dot\l_3=\dot\l_4=\dot6.$ For simplicity, we write
 $B=B_{6,6}^{6,6}$ and $\dot B=B^{\dot 6,\dot 6}_{\dot 6,\dot 6}.$

\bl{case2}
For $\dot\l_1=\dot\l_2=\dot\l_3=\dot\l_4=\dot6,$ the dotted equation (\ref{2}) holds by choosing
$$c_{\dot6,\dot\l_1}^{\dot \l_2}=c_{6,\l_1}^{\l_2}(T_{6,\l_1}^{\l_2})^{-1}$$
where $T_{6,\l_1}^{\l_2}$ are as in Lemma \ref{cross sym}.
\el
\proof
Lemma \ref{cross sym} asserts
$$B_{\mu,\gamma}=(T^6_{6,\mu})^{-1}(T^\mu_{6,6})^{-1}B_{\gamma,\mu}T^6_{6,\gamma}T^{\gamma}_{6,6}.$$
We also know $\dot B=(B^t)^{-1}.$ In this case, we rewrite equation (\ref{2}) as:
\begin{equation}\label{2*}
\Sigma_\mu c^6_{6,\mu}c^\mu_{6,6}B_{\mu,\gamma}=
\Sigma_{\mu} c^6_{6,\mu}c^\mu_{6,6}(T^6_{6,\mu})^{-1}(T^\mu_{6,6})^{-1}B_{\gamma,\mu}T^6_{6,\gamma}T^{\gamma}_{6,6}
=c^6_{6,\gamma}c^\gamma_{6,6},
\end{equation}
for
 $\mu, \gamma \in \{0,6\},$
i.e.
$$\Sigma_{\mu} c^6_{6,\mu}c^\mu_{6,6}(T^6_{6,\mu})^{-1}(T^\mu_{6,6})^{-1}B_{\gamma,\mu}
=c^6_{6,\gamma}c^\gamma_{6,6}(T^6_{6,\gamma})^{-1}(T^{\gamma}_{6,6})^{-1}.$$

Since $\dot B=(B^t)^{-1},$ obviously, the above equation implies
$$\Sigma_{\mu} c^6_{6,\mu}c^\mu_{6,6}(T^6_{6,\mu})^{-1}(T^\mu_{6,6})^{-1}\dot B_{\mu,\gamma}=c^6_{6,\gamma}c^\gamma_{6,6}(T^6_{6,\gamma})^{-1}(T^{\gamma}_{6,6})^{-1}.$$
In order to have dotted equation (\ref{2}), one only needs to take
$$c_{\dot6,\dot\l_1}^{\dot \l_2}=c_{6,\l_1}^{\l_2}(T_{6,\l_1}^{\l_2})^{-1}.$$
To ensure the skew-symmetry property of a vertex operator algebra, we need a normalization so that
$J^{\dot 6}(u,z)=\Y_{\dot6,\dot0}^{\dot6}+\frac{c_{\dot6,\dot6}^{\dot6}}{c_{\dot 6,\dot 0}^{\dot6}}\Y_{\dot6,\dot6}^{\dot6}+\frac{c_{\dot6,\dot6}^{\dot0}}{ c_{\dot 6,\dot0}^{\dot 6}}\Y_{\dot6,\dot6}^{\dot0}.$
\endproof

\begin{rmk}\label{voa nondeg}
One can see clearly that in the construction of $L_{\frak{sl}(10)}(2,0)^e,$ the only nontrivial braiding
matrix we use is $B^{6,6}_{6,6}.$ That is, we only use the fact that $D_{6,6}^\mu\neq 0,$ and $D_{6,\mu}^6\neq 0,$ for $\mu=0,2,4,6,8$ in the arguments.  As we mentioned already that a proof of the fact
using the language of \voa  is given in $\S$\ref{appendix}.  Similar calculation could be
done for $L_{\frak{sl}(28)}(2,0)^e.$
\end{rmk}

\begin{rmk}
By choosing
$\dot J^{\dot6}=\Y_{\dot6,\dot0}^{\dot6}
+a\frac{c_{\dot6,\dot6}^{\dot6}}{c_{\dot 6,\dot0}^{\dot6}}\Y_{\dot6,\dot6}^{\dot6}
+a^2\frac{c_{\dot6,\dot6}^{\dot0}}{ c_{\dot 6,\dot0}^{\dot 6}}\Y_{\dot6,\dot6}^{\dot0}$
for any $a\in \C^*$, the dotted equation (\ref{2}) is still satisfied.  We will first prove there
is a \voa structure on $L_{\frak{sl}(10)}(2,0)^e$ and then  in the next section prove that different choices of
the vertex operators actually give the same \voa structure.
\end{rmk}

\begin{lem}\label{rationalsol}
$E\langle \dot u_4, \dot J^{\dot{\l_3}}(\dot u_3,w)\dot J^{\dot{\l_2}}(\dot u_2,z)\dot u_1\rangle $
is rational, symmetric for all higher descendants.
\end{lem}

\proof
Symmetric property follows from Lemma \ref{rational} without the assumption that $\dot u_i$ is primary for $i=1,2,3,4.$

Let $\{x_{i,j}|,i,j=1,\cdots , 10,i\neq j\}\cup \{h_i|1\leq i\leq 9\}$ be the standard basis of $\frak{sl}(10)$
as in \cite{H-Lie}. Fix $\dot u_i$ i=1,2,3,4 as in Lemma \ref{rational}.
Assume that  $n=0$ and $i>j,$ or $n>0$ for any $i,j.$ Using equation (\ref{comm}), we have
\begin{equation}
\begin{aligned}
&E\langle x_{i,j}(-n)\dot u_4, \dot J^{\dot{\l_3}}(\dot u_3,w)\dot J^{\dot{\l_2}}(\dot u_2,z)\dot u_1\rangle\\
=&-E\langle\dot u_4,  x_{i,j}(n)\dot J^{\dot{\l_3}}(\dot u_3,w)\dot J^{\dot{\l_2}}(\dot u_2,z)\dot u_1\rangle\\
=&- E\langle\dot u_4,  \dot J^{\dot{\l_3}}(x_{i,j}(0)\dot u_3,w)\dot J^{\dot{\l_2}}(\dot u_2,z)\dot u_1\rangle
-E\langle\dot u_4,  \dot J^{\dot{\l_3}}(\dot u_3,w) x_{i,j}(n)\dot J^{\dot{\l_2}}(\dot u_2,z)\dot u_1\rangle\\
=&-E\langle\dot u_4,  \dot J^{\dot{\l_3}}(x_{i,j}(0)\dot u_3,w)\dot J^{\dot{\l_2}}(\dot u_2,z)\dot u_1\rangle
-E\langle\dot u_4,  \dot J^{\dot{\l_3}}(\dot u_3,w) \dot J^{\dot{\l_2}}(x_{i,j}(0)\dot u_2,z)\dot u_1\rangle \\
&~~~~~~~~~~~~~~~~-E\langle\dot u_4,  \dot J^{\dot{\l_3}}(\dot u_3,w) \dot J^{\dot{\l_2}}\dot u_2,z)(x_{i,j}(n)\dot u_1\rangle\\
=&-E\langle\dot u_4,  \dot J^{\dot{\l_3}}(x_{i,j}(0)\dot u_3,w)\dot J^{\dot{\l_2}}(\dot u_2,z)\dot u_1\rangle- E\langle\dot u_4,
 \dot J^{\dot{\l_3}}(\dot u_3,w) \dot J^{\dot{\l_2}}(x_{i,j}(0)\dot u_2,z)\dot u_1\rangle,
\end{aligned}
\end{equation}
which is a rational function by Lemmas \ref{case1} and \ref{case2}.

Using a similar argument, one easily gets
$$\langle \dot u_4, \dot J^{\dot{\l_3}}(\dot u_3,w)\dot J^{\dot{\l_2}}(\dot u_2,z) x_{i,j}(-n)\dot u_1\rangle$$
is also rational for either $n=0$ and $i<j$, or $n>0$.
Thus $\langle \dot u_4, \dot J^{\dot{\l_3}}(\dot u_3,w)\dot J^{\dot{\l_2}}(\dot u_2,z)\dot u_1\rangle$
 is rational for any $\dot{u}_4, \dot{u}_1\in L_{\frak{sl}(10)}(2,0)^e$ and primary elements $\dot u_2, $ $\dot u_3.$
 Together with the symmetric property this implies for any two primary elements $\dot u_2,\dot u_3$,
 there exists some $k\geq 0,$ such that
$$(w-z)^k [\dot J^{\dot{\l_3}}(\dot u_3,w),\dot J^{\dot{\l_2}}(\dot u_2,z)]=0$$
which is the locality condition.
Since the $L_{\frak{sl}(10)}(2,\dot{0})$-module  $L_{\frak{sl}(10)}(2,\dot{0})\oplus L_{\frak{sl}(10)}(2,\L_3+\L_7)$
 is generated by the primary elements,  it follows from Lemma \ref{Dong} that
for any $\dot u,\dot v \in L_{\frak{sl}(10)}(2,\dot{0})\oplus L_{\frak{sl}(10)}(2,\L_3+\L_7)$ there exists $r\geq 1$ such that
$$(z_1-z_2)^r [\dot J^{\dot{\l_3}}(\dot u,z_1),\dot J^{\dot{\l_2}}(\dot v,z_2)]=0,$$
as expected.
\endproof

Theorem \ref{main} now follows from Lemma \ref{rationalsol}.

\begin{rmk}
There have been some work in literature on  finding rational solutions of the
KZ-equations (cf. \cite{RST} and references therein) in special cases. Lemma \ref{rationalsol}
indicates that there is a rational solution of the KZ-equation for $SU(10)_2.$
\end{rmk}

\section{Uniqueness of $L_{\frak{sl}(10)}(2,0)^e$}\label{uni}
 We discuss the
uniqueness of the \voa structure on $L_{\frak{sl}(10)}(2,0)^e$ in this section. We first prove the uniqueness
of $L_{\frak{sl}(2)}(10,0)^e=L_{\frak{sl}(2)}(10,0)\oplus L_{\frak{sl}(2)}(10,6).$

\bl{B_2} Assume that $L_{\frak{sl}(2)}(10,0)^e=L_{\frak{sl}(2)}(10,0)\oplus L_{\frak{sl}(2)}(10,6)$ is a simple
vertex operator algebra which is an extension of $L_{\frak{sl}(2)}(10,0).$
Let $Y$ denote the vertex operator on the vertex operator algebra $L_{\frak{sl}(2)}(10,0)^e.$ As in $\S$\ref{key},
 for $u\in L_{\frak{sl}(2)}(10,6),$ set $$Y(u,z)= J^6(u,z)=\sum_{\l_1,\l_2\in\{0,6\}} d_{6,\l_1}^{\l_2} \mathcal{Y}_{6,\l_1}^{\l_2}(u,z),$$
then $$d_{6,\l_1}^{\l_2} \neq 0 \hbox{ if } \V_{6,\l_1}^{\l_2}\neq 0. $$
\el
\proof It is clear $d_{6,0}^6=c_{6,0}^6\neq 0$ by the skew symmetry. Note that $L_{\frak{sl}(2)}(10,0)^e$
has a non-degenerate, symmetric, invariant bilinear form $(\cdot,\cdot)$ (cf. \cite{Li0}). This implies that $d_{6,6}^0\ne 0.$

We now prove that $d_{6,6}^6\ne 0.$  Assume that $d_{6,6}^6= 0.$ It is known that the weight one subspace
$\frak{sl}(2)+L(6)$ of $L_{\frak{sl}(2)}(10,0)^e$ is a Lie algebra denoted by $\frak g$ where $L(6)$ is the irreducible
module with highest weight $6$ for $\frak{sl}(2).$ If $d_{6,6}^6= 0,$ then
$[L(6),L(6)]= \frak{sl}(2)$ or $0.$  If $[L(6),L(6)]= \frak{sl}(2),$ then $\frak g$ is a simple Lie algebra.
 According to the classification of finite dimensional simple Lie algebras,  the only possibility for
 $\frak g$ is $B_2,$ thus $L_{\frak{sl}(2)}(10,0)^e=L_{B_2}(1,0),$ a contradiction.

If $[L(6),L(6)]=0,$ then $L(6)$ generates a Heisenberg vertex operator algebra $U$ with central charge
$7.$ The character of the Heisenberg vertex operator algebra $U$ is
$$\ch_qU=\frac{q^{-\frac{5}{48}}}{\prod_{n\geq 1}(1-q^n)^7}$$
here $5/2$ is the central charge of $L_{\frak{sl}(2)}(10,0)$.
By applying Lemma \ref{polygrow} to $L_{\frak{sl}(2)}(10,0)+L_{\frak{sl}(2)}(10,6)$ as a $L_{\frak{sl}(2)}(10,0)$-module, we immediately get that the coefficients of
$$\eta(q)^{5/2}\ch_q (L_{\frak{sl}(2)}(10,0)+L_{\frak{sl}(2)}(10,6))$$
 satisfy the polynomial growth condition. But the coefficients of
$$\eta(q)^{5/2}\ch_qU=\frac{1}{\prod_{n\geq 1}(1-q^n)^{9/2}}$$
has exponential growth, a contradiction. The proof is complete.
\endproof

\begin{rmk}\label{Bnontri}
By Lemma \ref{B_2} and
\begin{eqnarray*}
& &  d^6_{6,0}d^0_{6,6}B_{0,0}+d^6_{6,6}d^6_{6,6}B_{6,0}=d^6_{6,0}d^0_{6,6}\\
& &  d^6_{6,0}d^0_{6,6}B_{0,6}+d^6_{6,6}d^6_{6,6}B_{6,6}=d^6_{6,6}d^6_{6,6}
\end{eqnarray*}
 which is an expansion of equation (\ref{2*}) with $c_{\lambda,\mu}^{\gamma}$ replaced by  $d_{\lambda,\mu}^{\gamma},$
  we see that the only option of $J^6$ that gives a vertex operator algebra structure on
the space $L_{\frak{sl}(2)}(10,0)^e$ is
$$J^6=c_{6,0}^6\Y_{6,0}^6+a^2(c_{6,6}^0\Y_{6,6}^0)+a(
c_{6,6}^6\Y_{6,6,}^6), \hbox{ for any } a\in \C^*.$$
\end{rmk}

The following Theorem will help us to determine the uniqueness of the \voa structure on both $L_{\frak{sl}(2)}(10,0)^e$ and
$L_{\frak{sl}(10)}(2,0)^e.$

\begin{thm}\label{unique}
Let $(V,Y,\1,\omega)$ be a \voa with a linear isomorphism $g$
 which preserves $\1$ and $\omega$. Set $Y^g(u,z)=g^{-1}Y(gu,z)g$ for any $u\in V.$
 Then $(V,Y^g,\1, \omega)$ is a \voa isomorphic to $(V,Y,\1,\omega).$

\end{thm}

\proof We first check the \voa axioms for $(V,Y^g,\1,\omega)$. Since $g \1=\1$ and
$g\omega =\omega , $ we only need to check the creativity, derivation
property and the commutativity.

1) Creativity: For $u\in V$
$$\lim_{z\to 0}Y^g(u,z)\1 =\lim_{z\to 0} g^{-1}Y(gu,z)g\1=g^{-1}\lim_{z\to 0} Y(gu,z)\1=g^{-1}gu=u .$$

2) Derivation property: Let $Y^g(\omega,z)=\sum_{n\in\Z}L^g(n)z^{-n-2}.$ Then
\begin{eqnarray*}
& & [L^g(n),Y^g(u,z)]=[g^{-1}L(-1)g,g^{-1}Y(gu,z)g]\\
& &\ \ \ \ \ =g^{-1}[L(-1),Y(gu,z)]g\\
& &\ \ \ \ \ =g^{-1}\frac{d}{dz}Y(gu,z)g\\
& &\ \ \ \ \ =\frac{d}{dz}g^{-1}Y(gu,z)g\\
& &\ \ \ \ \ =\frac{d}{dz}Y^g(u,z).
\end{eqnarray*}

3) Commutativity: For any $u,z \in V, $ by commutativity of $(V,Y,\1,\omega),$ there exists $n\in \Z$
such that
$$(z_1-z_2)^n[Y(gu,z_1),Y(gv,z_2)]=0.$$
This implies that  $$(z_1-z_2)^n[Y^g(u,z_1),Y^g(v,z_2)]=0.$$
Thus  $(V,Y^g,\1,\omega)$ is a vertex operator algebra. It is clear that the linear map
 $g:V\to V$ gives a \voa isomorphism from $(V,Y^g,\1, \omega)$ to $(V,Y,\1, \omega).$
\endproof

\begin{cor}
 The simple vertex operator algebra structure on $L_{\frak{sl}(2)}(10,0)^e$ is unique.
 \end{cor}
\proof
Let $V=L_{\frak{sl}(2)}(10,0)+L_{\frak{sl}(2)}(10,6)$ be a simple vertex operator algebra
which is an extension of $L_{\frak{sl}(2)}(10,0).$ Then
$J^6=c_{6,0}^6\Y_{6,0}^6+a^2(c_{6,6}^0\Y_{6,6}^0)+a(c_{6,6}^6\Y_{6,6,}^6),$
for some $a\in \C^*.$ Note that $L_{B_2}(1,0)$ and $V$ are isomorphic $L_{\frak{sl}(2)}(10,0)$-modules.
Let $g: L_{B_2}(1,0)\to L_{B_2}(1,0)$
 be the linear map such that $g|_{L_{\frak{sl}(2)}(10,0)}=1$ and $g|_{L_{\frak{sl}(2)}(10,6)}=a.$ Then
 $V$ and $L_{B_2}(1,0)^g$ are isomorphic by noting that for $u\in L_{\frak{sl}(2)}(10,6),$
$$ Y^g(u,z)=c_{6,0}^6\Y_{6,0}^{6}(u,z)+a c_{6,6}^{6}\Y_{6,6}^{6}(u,z)+a^2c_{6,6}^{0}\Y_{6,6}^{0}(u,z)$$
(see Remark \ref{Bnontri}). Thus, $V$ and $L_{B_2}(1,0)$ are isomorphic by Theorem \ref{unique}.
\endproof

The following corollary follows from a similar argument.

\begin{cor}\label{Y^g}
The simple vertex operator algebra structure on $L_{\frak{sl}(10)}(2,0)^e$ is unique.
\end{cor}

%----------------------------------------------------------------------------------
%----------------------------------------------------------------------------------

\section{Mirror extension of $L_{\frak{sl}(28)}(2,0)$}\label{28}

We give another example of mirror extension which is the extension of
$L_{\frak{sl}(28)}(2,0)$ in this section. Although this example
is more complicated than the example given in Section \ref{key}, the ideals and the
methods are similar.

Consider the conformal inclusion $SU(2)_{28}\subset (G_2)_1$
 (see \cite{CIZ, GNO}) and the level-rank duality $SU(2)_{28}\times SU(28)_2\subset SU(56)_1.$
Due to \cite{CIZ, DLN}, one knows that the \voa
\be{G_2}
L_{G_2}(1,0)=L_{\frak{sl}(2)}(28,0)\oplus L_{\frak{sl}(2)}(28,10)\oplus L_{\frak{sl}(2)}(28,18)\oplus L_{\frak{sl}(2)}(28,28),
\ee
 where $L_{\frak{sl}(2)}(28,0)\subset L_{G_2}(1,0)$ is a conformal embedding with central charge $\frac{14}{5}.$

  The decomposition of
 $L_{\frak{sl}(56)}(1,0)$ under $L_{\frak{sl}(2)}(28,0) \otimes L_{\frak{sl}(28)}(2,0)$ is given in Lemma \ref{duality} as follow
\be{56}
\begin{aligned}
L_{\frak{sl}(56)}(1,0)=\bigoplus_{ a=0, ~a~ \hbox{even,}}^{28} a\times \dot a,
\end{aligned}
\ee
where $a=L_{\frak{sl}(2)}(28,a)$ and
$\dot a=L_{\frak{sl}(28)}(2,\Lambda_{\frac{a}{2}}+\Lambda_{28-\frac{a}{2}})$, here $\Lambda_i$
 is the fundamental weight of ${\frak{sl}}(28),$ and we use $\Lambda_{28}=\Lambda_0$
  (sometimes 0 by abusing of notations) to denote the trivial representation of $\frak{sl}(28).$

By the Mirror Extension Conjecture and equations (\ref{G_2}) and (\ref{56}), it is expected that
there is a vertex operator algebra structure on
$$V=L_{\frak{sl}(28)}(2,0)^e=L_{\frak{sl}(28)}(2,0)+L_{\frak{sl}(28)}(2,\Lambda_5+\Lambda_{23})
+L_{\frak{sl}(28)}(2,\Lambda_9+\Lambda_{19})+L_{\frak{sl}(28)}(2,2\Lambda_{14})$$
with central charge $\frac{261}{5}.$  Note that the vertex operator algebra structure on $V$
cannot be obtained from the framed vertex operator algebras as $\frac{261}{5}$ is not a half integer.
For conveniece, we use $V=\dot 0+\dot {10}+\dot{18}+\dot{28}.$

\begin{thm}
There is a \voa structure on $V=L_{\frak{sl}(28)}(2,0)^e.$
\end{thm}

\proof
Again, we only need to define the vertex operator $\dot Y$ on $V$ satisfying locality.

We introduce some notations first. We use $\tilde Y,$ $Y$ and $\dot Y$ to denote the
vertex operators on $L_{\frak{sl}(56)}(1,0),$ $L_{G_2}(1,0)$ and $V$ respectively.
As in $\S$\ref{key} ,  for $u_1\in L_{\frak{sl}(2)}(28,\l),$ $u_2\in L_{\frak{sl}(28)}(2,\dot{\l}),$
\be{VO SU(56)}
\tilde Y(u_1\otimes u_2,z)=\Sigma_{\l_1, \l_2}D^{\l_2}_{\l,\l_1}\cdot \mathcal{\Y}^{\l_2}_{\l,\l_1}(u_1,z)
\otimes \mathcal{Y}^{\dot{\l_2}}_{\dot{\l},\dot{\l_1}}(u_2,z),
\ee
\be{VO g2}
Y(u_1,z)=\sum_{\l_1,\l_2} c_{\l,\l_1}^{\l_2} \Y_{\l,\l_1}^{\l_2}(u_1,z)
\ee
where  $\Y^{\l_2}_{\l,\l_1}\in\V^{\l_2}_{\l,\l_1},$ $\Y^{\dot \l_2}_{\dot \l,\dot \l_1}\in\V^{\dot \l_2}_{\dot\l,\dot\l_1}.$
%$\V^{\l_1}_{\l,\l_3}$ is the space of intertwining operators of type $\left (
%\begin{array}{c}
%\ \ \l_3 \ \\
%\l \ \  \l_1
%\end{array}
 %\right ),$ and $D_{\l,\l_3}^{\l_1}\in \C .$
 As in equation (\ref{VO SU(20)}),  we can assume that $D_{\l,\l_1}^{\l_2}=1$ whenever $\V_{\l,\l_1}^{\l_2}\neq 0$
 for $\l, \l_1, \l_2 \in \{0,10,18,28\}$
 by suitably choosing the basis of intertwining operators
 $\Y_{\l,\l_1}^{\l_2}$ and $\dot \Y_{\dot\l,\dot\l_1}^{\dot\l_2}.$  This can also be  calculated in the frame
  work of \voa  similarly as in $\S$\ref{appendix}.

We now determine the vertex operator $\dot Y$ on $V.$
For $u\in L_{\frak{sl}(28)}(2,\dot \l),$ as
 in $\S$\ref{key},  we write
$$
\dot Y(u,z)=\dot J^{\dot \l}(u,z) =\sum_{\dot\l_1,\dot \l_2}\dot c_{\dot\l,\dot\l_1}^{\dot\l_2}\Y_{\dot\l,\dot\l_1}^{\dot\l_2}.$$
We need the coefficients satisfy an equation similar to the dotted equation (\ref{2}).

For $u\in L_{\frak{sl}(28)}(2,0),$ the choice of $\dot Y(u,z)$ is obviously. For
 $u\in L_{\frak{sl}(28)}(2,\Lambda_5+\Lambda_{23}), $ same as in Lemma \ref{case2}, we can take
$$\dot c_{\dot{10},\dot\l_1}^{\dot\l_2}=c_{10,\l_1}^{\l_2}(T_{10,\l_1}^{\l_2})^{-1}$$
 to  guarantee locality.

Since $L_{\frak{sl}(28)}(2,2\Lambda_{14})$ is a simple current, the braiding matrix
$B^{\dot{28},\dot\l_2}_{\dot \l_4,\dot\l_1} $ is just a number.
Equations similar to (\ref{2}) and (\ref{eq2}) imply that
$$B^{{28},\l_2}_{ \l_4,\l_1}=(B^{\l_2,{28}}_{\l_4,\l_1})^{-1}=B^{\dot{\l}_2,\dot{28}}_{\dot \l_4,\dot\l_1} .$$
 %and
% $$B^{\dot{28},\dot\l_2}_{\dot \l_4,\dot\l_1}=B^{\l_2,{28}}_{\l_4,\l_1}$$ for $\l_2\in \{0,10,28\}.$
Using
$$c_{\l_2,\mu}^{\l_4}c_{28,\l_1}^{\mu}B^{\l_2,28}_{\l_4,\l_1}=c_{28,\gamma}^{\l_4}c_{\l_2,\l_1}^{\gamma}$$
and equation (\ref{cross sym}), we get
\be{}
 c_{\l_2,\mu}^{\l_4}c_{28,\l_1}^{\mu} (T_{\l_2,\mu}^{\l_4})^{-1}(T_{28,\l_1}^{\mu})^{-1}B^{28,\l_2}_{\l_4,\l_1}=
 c_{28,\gamma}^{\l_4}c_{\l_2,\l_1}^{\gamma}(T_{28,\gamma}^{\l_4})^{-1}(T_{\l_2,\l_1}^{\gamma})^{-1},
\ee
or equivalently,
$$ c_{\l_2,\mu}^{\l_4}c_{28,\l_1}^{\mu} (T_{\l_2,\mu}^{\l_4})^{-1}(T_{28,\l_1}^{\mu})^{-1}B^{\dot{\l}_2,\dot{28}}_{\dot\l_4,\dot\l_1}=
c_{28,\gamma}^{\l_4}c_{\l_2,\l_1}^{\gamma}(T_{28,\gamma}^{\l_4})^{-1}(T_{\l_2,\l_1}^{\gamma})^{-1}.$$
As long as we choose $c_{\dot{28},\dot\l_3}^{\dot\l_1}=c_{28,\l_3}^{\l_1}(T_{28,\l_3}^{\l_1})^{-1},$
we get $\dot J^{\dot {28}},$ $\dot J^{\dot 0}$ and $\dot J^{\dot {10}},$ which are pairwise mutually local.

We have defined $\dot Y (u,z)$ for
$u\in L_{\frak{sl}(28)}(2,0)\oplus L_{\frak{sl}(28)}(2,\Lambda_{5}
+\Lambda_{23})\oplus L_{\frak{sl}(28)}(2,2\Lambda_{14}).$
It remains to define
$\dot Y (u,z)$ for $u\in L_{\frak{sl}(28)}(2,\Lambda_{9}+\Lambda_{19}).$
Since the fusion  product
$$L_{\frak{sl}(28)}(2,\Lambda_{5}+\Lambda_{23})\boxtimes
L_{\frak{sl}(28)}(2,2\Lambda_{14})=L_{\frak{sl}(28)}(2,\Lambda_{9}+\Lambda_{19}),$$
 any
$u\in L_{\frak{sl}(28)}(2,\Lambda_{9}+\Lambda_{19})$
can be expressed
as $u=\sum_i(a^i)_{m_i}b^i,$ for some
$a^i\in L_{\frak{sl}(28)}(2,\Lambda_{5}+\Lambda_{23}),$
and $b^i\in L_{\frak{sl}(28)}(2,2\Lambda_{14}).$
 Thus for $u\in \dot{18}$ one can define
$$\dot J^{\dot{18}}(u,z)=\sum_i (\dot J^{\dot{10}}(a^i,z))_{m_i}\dot J^{\dot{28}}(b^i,z)$$
(see  equation (\ref{OPE})). Lemma \ref{Dong} ensures locality of $\dot Y$ defined on $V.$
Thus $(V,\dot Y,\1,\omega)$ gives a \voa structure on $V.$
\endproof

For the uniqueness of the structure,  it is quite similar to the proof of uniqueness of
$L_{\frak{sl}(10)}(2,0)^e$ by viewing $V$ as an extension of the rational and $C_2$-cofinite \vosa
$$U=L_{\frak{sl}(10)}(2,0)\oplus L_{\frak{sl}(10)}(2,2\Lambda_{14}).$$ Then
$M=L_{\frak{sl}(10)}(2,\Lambda_{5}+\Lambda_{23})\oplus L_{\frak{sl}(10)}(2,\Lambda_{9}+\Lambda_{19})$ is an irreducible
 $U$-module.  Since the  structure of $U$ is unique and $M$ also
 has a unique $U$-module structure,  we derive $V=U+M$ has a unique \voa structure as in  Corollary \ref{Y^g}.

\section{Comments on general case}\label{comments}

The general idea presented is in principle applicable to higher rank case. However,
there are a number of technical problems which are not resolved in the literature. For example,
for $SU(n)_k$ with $n\geq 3,$ it is not clear if  the braiding matrices coming
from solutions of KZ equation are similar to unitary matrices.

From categorical point of view, Theorem 3.8 in \cite{Xu-m} can be seen as a statement
about existence of commutative Frobenius  algebras from given ones. For example, in
the case of the key example, theorem 3.8 in \cite{Xu-m} says that $\dot{0}+\dot{6}$ is a
commutative Frobenius algebra in the unitary tensor category associated with $SU(10)_2.$
According to  \cite{HK}, this is equivalent to the existence of
local extensions of $SU(10)_2$ with spectrum $\dot{0}+\dot{6}.$ However, to apply
\cite{HK} one must show that the unitary tensor  category associated with
$SU(10)_2$ from the operator algebra framework is the same as that coming from the
theory of vertex operator algebra. In the case of $SU(10)_2$ one can presumably use
the cohomology vanishing argument in \cite{Yasu-L}. But this is not entirely clear, since
the braiding matrix in operator algebra is automatically unitary, and we do not even know
if for $SU(10)_2$ case, the braiding matrix from solutions of KZ equationis is similar to a unitary matrix.

The rationality of both $L_{\frak{sl}(10)}(2,0)^e$ and $L_{\frak{sl}(28)}(2,0)^e$ have not been
 established in this paper. They are completely rational in operator algebra framework, and
 in fact all its irreducible representations of $L_{\frak{sl}(10)}(2,0)^e$  are listed on P. 96 of
 \cite{Xu-am}, where their irreducible representations are used with simple current extensions
 to construct holomorphic $c=24$ nets which
 corresponds to number 40 in Schelleken's list (\cite{S}).  It is worthy to point out that
 some holomorphic $c=24$ \voas including number 40 in
 Schelleken's list are constructed  in \cite{L, LS} by using framed \voas. We are informed recently
 that the \voa $L_{\frak{sl}(10)}(2,0)^e $ can be realized as a coset vertex operator algebra in the framed holomorphic \voa
 \cite{L} corresponding to number 40 in Schelleken's list. We plan to investigate the connection of $L_{\frak{sl}(10)}(2,0)^e $ with the framed vertex operator algebras further.

%----------------------------------------------------------------------------------
%----------------------------------------------------------------------------------
\section{Appendix}\label{appendix}
We prove the claim made in Section 4 (see Remark \ref{voa nondeg}) in the Appendix. That is, $D_{6,6}^{\mu}\ne 0,$ $D_{6,\mu}^6\ne 0$ for $\mu=0,2,4,6,8$ where $D_{\l,\l_1}^{\l_2}$ are defined in (\ref{VO SU(20)}).

\subsection{Decomposition of $L_{\frak{sl}(20)}(1,0)$ as $L_{\frak{sl}(2)}(10,0)\otimes L_{\frak{sl}(10)}(2,0)$-module}

Let $V_L$ denote the lattice \voa   associated to the root lattice of $\frak{sl}(20)$ $(V_L\cong L_{\frak{sl}(20)}(1,0)$).
We use $\{\e_1,\cdots , \e_{20}\}$ to denote the standard orthonormal basis of
$\R^{20}$ with the usual inner product.  Then $L=\Sigma_{i\neq j,i,j=1 }^{20}\Z (\e_i-\e_j).$
 Set $\a_i=\e_i-\e_{i+1},$ $1\leq i\leq 19.$
The lattice \voa $V_L=M(1)\otimes \C^\varepsilon[L],$ where $\varepsilon: L\times L \to \langle \pm 1\rangle$
is a 2-cocycle s.t. $\varepsilon (\a, \b)\varepsilon(\b, \a)=(-1)^{\a,\b},$ $\forall~ \a, \b \in L.$
Set
$$x_{i,j}=E_{i,j}+E_{10+i,10+j},~h_{i,j}=E_{ii}-E_{jj}+E_{10+i,10+i}-E_{10+j,10+j},$$
for $1\leq i,j\leq10,$ $ i \neq j,$ where $E_{i,j}$ denotes the $20\times 20$ matrix with $1$ in the $(i,j)$-entry
and $0$ elsewhere.  The bilinear form is defined as
$(A, B)=\frac{1}{2}tr (AB),$ thus $(h_{i,j},h_{i,j})=4.$
Then
$$\Sigma _{i\neq j}\C x_{i,j}+\Sigma_{i\neq j}\C h_{i,j}\cong \frak{sl}(10)\subset \frak{sl}(20).$$
This gives a vertex operator algebra embedding  $L_{\frak{sl}(10)}(2,0)\subset V_L.$
 We use
 $\b_{i}=E_{i,i}-E_{i+1,i+1}+E_{10+i,10+i}-E_{10+i+1,10+i+1},$ $1\leq i\leq 9,$ to
 denote a basis of the Cartan algebra of the sub Lie algebra $\frak{sl}(10).$

Set
\begin{equation}\label{sl_2}
\begin{split}
&e=E_{1,11}+E_{2,12}+\cdots +E_{10,20};\\
&f=E_{11,1}+E_{12,2}+\cdot +E_{20,10};\\
&h=E_{1,1}+\cdots+E_{10,10}-E_{11,11}-\cdots-E_{20,20},
\end{split}
\end{equation}
$\C e+\C f +\C h\cong \frak{sl}(2)\subset \frak{sl}(20).$
Notice that $(h,h)=20,$ one immediately see that the \voa $V_L$ has a \vosa
$L_{\frak{sl}(2)}(10,0)$ associated to the inclusion $\frak{sl}(2)\subset \frak{sl}(20).$

It is easy to check that  $L_{\frak{sl}(10)}(2,0)\subset L_{\frak{sl}(2)}(10,0)^c$ (the commutant of $L_{\frak{sl}(2)}(10,0)$ in $V_L$).
Thus on \voa  level, we have an inclusion
$L_{\frak{sl}(10)}(2,0)\otimes L_{\frak{sl}(2)}(10,0)\subset  V_L.$
By considering the central charge,
\begin{equation*}
\begin{split}
&c_{L_{\frak{sl}(2)}(10,0)}=\frac{30}{10+2}=5/2,\\
&c_{L_{\frak{sl}(10)}(2,0)}=\frac{2\times 99}{10+2}=33/2,\\
&c_{V_L}=19=33/2+5/2,
\end{split}
\end{equation*}
we get the inclusion is actually a conformal inclusion.

Equation (\ref{eqsl20}) gives  the decomposition of $V_L$ as
a $ L_{\frak{sl}(2)}(10,0) \otimes L_{\frak{sl}(10)}(2,0)$ module :
\begin{equation}\label{decompo}
\begin{split}
V_L&=  L_{\frak{sl}(2)}(10,0) \otimes L_{\frak{sl}(10)}(2,0)\\
&\oplus L_{\frak{sl}(2)}(10,2) \otimes L_{\frak{sl}(10)}(2,\Lambda_1+\Lambda_9)\\
&\oplus L_{\frak{sl}(2)}(10,4) \otimes L_{\frak{sl}(10)}(2,\Lambda_2+\Lambda_8)\\
&\oplus L_{\frak{sl}(2)}(10,6) \otimes L_{\frak{sl}(10)}(2,\Lambda_3+\Lambda_7)\\
&\oplus L_{\frak{sl}(2)}(10,8) \otimes L_{\frak{sl}(10)}(2,\Lambda_4+\Lambda_6)\\
&\oplus L_{\frak{sl}(2)}(10,10) \otimes L_{\frak{sl}(10)}(2,2\Lambda_5).
\end{split}
\end{equation}
We denote the decomposition as
$0\times \dot{0}+2\times \dot{2}+4\times \dot{4}+6\times \dot{6}+8\times \dot{8}+10\times \dot{10}.$ For short we set  $M^{\l}=\l\times\dot \l.$
Note that each $M^{\lambda}=\oplus_{n=0}^{\infty}M^{\l}(n)$ is an irreducible  highest weight module for the affine algebra $\hat{\frak{sl}(2)}\times \hat{\frak{sl}(10)}$-module and $M^{\l}(0)\cong L_{\frak{sl}(2)}(\l)\otimes L_{\frak{sl}(10)}(\dot\l)$ is an irreducible  $\frak{sl}(2)\times \frak{sl}(10)$-module where $L_{\frak{g}}(\lambda)$ is the irreducible highest weight module for a finite dimensional simple Lie algebra $\frak g$ with highest weight
$\lambda.$

We now determine the highest (lowest, resp.)  weight vectors of $M^{\lambda}$ which are the highest (lowest, resp.) weight vectors of $\frak{sl}(2)\times \frak{sl}(10)$-modules of $M^{\l}(0).$

Since $$(V_L)_1\cong \frak{sl}(20)\supset \frak{sl}(2)\oplus \frak{sl}(10),$$
 we also use $x_{i,j},$ $h_{i,j},$ $e, f$ and $h$ for
elements in $(V_L)_1.$ Set
\begin{equation}
\begin{split}
&v^1=e^{\e_1-\e_{20}}, \\
&v^2=e^{\e_1+\e_2-\e_{19}-\e_{20}}, \\
&v^3=e^{\e_1+\e_2+\e_3-\e_{18}-\e_{19}-\e_{20}},\\
&v^4=e^{\e_1+\e_2+\e_3+\e_4-e_{17}-\e_{18}-\e_{19}-\e_{20}},\\
&v^5=e^{\e_1+\e_2+\e_3+\e_4+\e_5-\e_{16}-\e_{17}-\e_{18}-\e_{19}-\e_{20}}.
\end{split}
\end{equation}
We claim that $v^i\in M^{2i}$ are the highest weight vectors. Since the proof is similar for all $i$ we
just demonstrate the proof for $i=3.$

It is easy to see
$$h_0v^3=6v^3, ~~(\b_7)_0v^3=(\b_3)_0v^3=v^3,$$ and $(\b_i)_0v^3=0$ for $i\neq 3,7.$
Since
\begin{equation*}
\begin{split}
Y(e,z)v^3&=\Sigma_{i=1}^{10}Y(e^{\e_i-\e_{10+i}},z)v^3\\
&=\Sigma_{i=1}^{10}E^-(\e_{10+i}-\e_i,z)E^
+(\e_{10+i}-\e_i,z)e_{\e_i-\e_{10+i}}z^{\e_i-\e_{10+i}}v^3\\
&=z\Sigma_{i=1}^{3}E^-(\e_{10+i}-\e_i,z)e_{\e_i-\e_{10+i}} v^3+\Sigma_{i=4}^{7}E^-(\e_{10+i}-\e_i,z)e_{\e_i-\e_{10+i}} v^3\\
&+z\Sigma_{i=8}^{10}E^-(\e_{10+i}-\e_i,z)e_{\e_i-\e_{10+i}} v^3
\end{split}
\end{equation*}
and
\begin{equation*}
\begin{split}
  ~~Y(f,z)v^3&=\Sigma_{i=1}^{10}Y(e^{\e_{10+i}-\e_i},z)v^3\\
&=\Sigma_{i=1}^{10}E^-(\e_i-\e_{10+i},z)E^
+(\e_i-\e_{10+i},z)e_{\e_{10+i}-\e_i}z^{\e_i-\e_{10+i}}v^3\\
&=z^{-1}\Sigma_{i=1}^{3}E^-(\e_i-\e_{10+i},z)e_{\e_{10+i}-\e_i} v^3+\Sigma_{i=4}^{7}E^-(\e_i-\e_{10+i},z)e_{\e_{10+i}-\e_i} v^3\\
&+z^{-1}\Sigma_{i=8}^{10}E^-(\e_i-\e_{10+i},z)e_{\e_{10+i}-\e_i} v^3,
\end{split}
\end{equation*}
we have
$$e_nv^3=0, ~ \forall~ n\geq 0, \hbox{ and } f_nv^3=0,~~\forall ~n>0.$$
Similarly,
$$(x_{i,i+1})_nv^3=0,\hbox{ if } n\geq 0, \hbox{ and }(x_{i+1,i})_nv^3=0, \hbox{ for } n>0.$$
Thus $v^3$ is a highest weight vector of $M^6.$

Similarly,
\begin{equation}
\begin{split}
&v_1=e^{-\e_1+\e_{20}},\\
&v_2=e^{-\e_1-\e_2+\e_{19}+\e_{20}}, \\
&v_3=e^{-\e_1-\e_2-\e_3+\e_{18}+\e_{19}+\e_{20}},\\
&v_4=e^{-\e_1-\e_2-\e_3-\e_4+e_{17}+\e_{18}+\e_{19}+\e_{20}},\\
&v_5=e^{-\e_1-\e_2-\e_3-\e_4-\e_5+\e_{16}+\e_{17}+\e_{18}+\e_{19}+\e_{20}}
\end{split}
\end{equation}
are the lowest weight vectors.

\subsection{Non-vanishing of $D_{6,\l_1}^{\l_2}$}
 We can now  prove that $D_{6,6}^\mu\neq 0,$ for $\mu=0,2,4,6,8$ and $D_{6,\gamma}^6\neq 0$
 for $\gamma =0,2,4,6,8$ (see Remark \ref{voa nondeg}). Recall equations (\ref{decompo}) and (\ref{VO SU(20)}). Then
 the vertex operator $Y(u,z)$ for $u\in M^{\l}$ on $V_L$ can be written as
 $$Y(u,z)=\sum_{\l_1, \l_2\in \{0,2,4,6,8,10\}}D_{\l,\l_1}^{\l_2}\Y^{\l_2}_{\l,\l_1}(u,z)
$$
where $\Y^{\l_2}_{\l,\l_1}$ is a fixed intertwining operator of type $\left (
\begin{array}{c}
\ \ M^{\l_2} \ \\
M^{\l} \ \  M^{\l_1}
\end{array}
 \right ).$

 \bl{simplify}
 If $D_{\l,\l_1}^{\l_2}\neq 0,$  then $D_{\l_1,\l}^{\l_2}\neq 0,$ $D_{\l,\l_2}^{\l_1}\neq 0$ and $D_{10-\l,10-\l_1}^{\l_2}\neq 0.$
 \el
 \proof  For $u\in M^\l,$ we denote $Y(u,z)$ by $\Y^\l(u,z).$
  Assume $D_{\l,\l_1}^{\l_2}\neq 0.$ Using the skew-symmetry and the fact that $V_L$ is self dual, we conclude
 $D_{\l_1,\l}^{\l_2}\neq 0,$ $D_{\l,\l_2}^{\l_1}\neq 0.$

 Since $D_{\l,\l_1}^{\l_2}\neq 0,$ there exist $u\in M^\l$ and $v \in M^{\l_1}$ such that
 $u_m v$ has a nonzero projection to $M^{\l_2}$ i.e. $$\la M^{\l_2}, \Y^\l(u,z)M^{\l_1}\ra\neq 0$$
 where we have used the fact that each $M^{\lambda}$ is a self dual $L_{\frak{sl}(2)}(10,0) \otimes L_{\frak{sl}(10)}(2,0)$-module.

 Fix nonzero homogeneous $b\in M^{10}.$ Then $M^0=\<a_nb|a\in M^{10}, n\in\Z\>.$ So we can find some
  $a\in M^{10}$ and $m\in\Z$
  such that
  $ a_m b= \1.$ We have
$$ \la M^{\l_2} , \Y^\l(u,z_1) \Y^{10}(a,z_2)\Y^{10}(b,z_3)M^{\l_1}\ra \ne 0.$$
 Using the associativity of vertex operators we see that
 $$\la M^{\l_2} ,\Y^{10-\l}(\Y^{\l}(u,z_1-z_2)a,z_2)\Y^{10}(b,z_3)M^{\l_1}\ra  \neq 0.$$
This implies $D_{10-\l,10-\l_1}^{\l_2}\neq 0.$
%\begin{equation*}
 %\begin{aligned}
 %u_m v &= u_m (a_nb)_{-1}  ~v\\
% &= u_m\sum_{i\geq 0}  {n \choose i}(-1)^ia_{n-i}b_{i-1}   ~v+(-1)^{n+1}u_m\sum_{i\geq 0}{n \choose i}(-1)^ib_{n-1-i}a_i~v\\
% &= u_m\sum_j (a^j)_{m_j}(b^j)_{n_j}~v,
%\end{aligned}
 %\end{equation*}
 %for some $a^j, b^j \in M^{10}.$
 %\begin{equation*}
 %\begin{aligned}
 %u_m v &= u_m\sum_j (a^j)_{m_j}(b^j)_{n_j}v\\
% &=\sum_j (u_m (a^j)_{m_j})(b^j)_{n_j}v\\
 %&=\sum_{j} [u_m, (a^j)_{m_j}](b^j)_{n_j}v+\sum_j (a^j)_{m_j}u_m(b^j)_{n_j}v
%\end{aligned}
% \end{equation*}
 \endproof

Thanks to the Lemma above, we only need to determine that $D_{6,6}^2,~ D_{6,6}^4, ~D_{6,6}^6$ and $D_{6,6}^8$ are nonzero (or $D_{6,6}^2, ~D_{4,4}^4,$ $D_{4,4}^6$ and $D_{4,4}^8$ are nonzero). There are several cases.
\bigskip

{\bf 1. $D_{6,6}^2\neq 0$:} We need to find $u, v \in M^6$ such that $u_m v  $ has a nonzero projection
to $M^2.$
 Take $f$ as in equation (\ref{sl_2}),  then
$$f_0v^1=e^{\e_{11}-\e_{20}}-e^{\e_{1}-\e_{10}}\in M^2.$$
 We want to show $e^{\e_{11}-\e_{20}}$ and $ e^{\e_{1}-\e_{10}}$ are in $M^6\cdot M^6,$
 where $M^i\cdot M^j= \la u_mv| u\in M^i, v\in M^j, m\in \Z\ra.$

Direct calculations give
\be{use}
(x_{1,7})_0v_3 =e^{-\e_2-\e_3-\e_7+\e_{18}+\e_{19}+\e_{20}}\in M^6,
\ee
 $$(x_{7,10})_0 ((x_{1,7})_0v_3)_4 v^{3}=C e^{\e_1-\e_{10}}\in M^6\cdot M^6,$$ for some $C \in \{\pm 1\}.$
Similarly, $e^{\e_{11}-\e_{20}}\in M^6\cdot M^6.$ Thus $f_0 v^1\in M^6\cdot M^6$ and $D_{6,6}^2\neq 0.$
 \bigskip

Before we deal with the other cases, we need the following lemma which is immediate using the proof of
$e^{\e_1-\e_{10}}\in M^6\cdot M^6.$
\bl{vectors}
Any element of type $
e^{\pm\e_{i_1}\pm\cdots \pm\e_{i_\l}-\mp\e_{10+j_1}\mp\cdots \mp\e_{10+j_\l} }$
lies in $M^{2\l},$ for $\l=0,1,2,3,4,$ %$ for $k+s=l+r=\l,$
where $i_k, j_l $ are distinct numbers in $\{1,\cdots, 10\}.$
\el

{\bf 2. $D_{4,4}^8\neq 0$:}
It is easy to see
$$a=e^{\e_3+\e_4-\e_{17}-\e_{18}}\in M^4$$ by considering
$(x_{3,1})_0(x_{4,2})_0(x_{7,9})_0(x_{8,10})_0 v^3$ (or Lemma \ref{vectors}).
Direct calculation gives
$$a_{-1} v^2=C_1 v^4\in M^8, $$
for some $C_1 \in \{\pm 1\},$
 which implies $D_{4,4}^8\neq 0.$

\bigskip

{\bf 3. $D_{4,4}^4\neq 0$: }
Notice that
\be{}
\begin{split}
 M^4\ni f_0f_0v^2+4v^2=&2(e^{\e_1+\e_2-\e_9-\e_{10}}+e^{\e_{11}+\e_{12}-\e_{19}-\e_{20}}+e^{\e_1-\e_9+\e_{12}-\e_{20}}\\
 &-e^{\e_1-\e_{10}+\e_{12}-\e_{19}}-e^{\e_2-\e_9+\e_{11}-\e_{20}}+e^{\e_2-\e_{10}+\e_{11}-\e_{20}}).
\end{split}
\ee

We need to show $e^{\e_1+\e_2-\e_9-\e_{10}},$ $e^{\e_{11}+\e_{12}-\e_{19}-\e_{20}},$
$e^{\e_1-\e_9+\e_{12}-\e_{20}},$ $e^{\e_1-\e_{10}+\e_{12}-\e_{19}},$
$e^{\e_2-\e_9+\e_{11}-\e_{20}},$ $e^{\e_2-\e_{10}+\e_{11}-\e_{20}}\in M^4\cdot M^4.$
We first prove that $e^{\e_1+\e_2-\e_3-\e_4}$ and $e^{\e_1-\e_3+\e_{18}-\e_{20}}$
 are in $M^4\cdot M^4.$ By Lemma \ref{vectors}, we know that
 $a=e^{-\e_2-\e_3+\e_{18}+\e_{19}}, b=e^{-\e_3-\e_4+\e_{19}+\e_{20}}$ lie in $M^4.$ Then
 \be{}
 \begin{aligned}
 &a_1v^2=C_3 e^{\e_1-\e_3+\e_{18}-\e_{20}}\in M^4\cdot M^4,\\
 & b_1 v^2= C_4 e^{\e_1+\e_2-\e_3-\e_4}\in M^4\cdot M^4.
 \end{aligned}
 \ee for some $C_3, C_4 \in \{\pm 1\}.$  Suitably choose some $(x_{i,j})_0$ acting on $e^{\e_1+\e_2-\e_3-\e_4}$
 we can get  $$e^{\e_1+\e_2-\e_9-\e_{10}}\in M^4\cdot M^4.$$ Similarly, $$e^{\e_{11}+\e_{12}-\e_{19}-\e_{20}}\in M^4\cdot M^4.$$
Suitable choosing $(x_{i,j})_0$ acting on $e^{\e_1-\e_3+\e_{18}-\e_{20}}$ asserts that
  $$e^{\e_1-\e_9+\e_{12}-\e_{20}},e^{\e_1-\e_{10}+\e_{12}-\e_{19}},
  e^{\e_2-\e_9+\e_{11}-\e_{20}}, e^{\e_2-\e_{10}+\e_{11}-\e_{20}}\in M^4\cdot M^4.$$
 This implies that $f_0f_0v^2+4v^2\in M^4\cdot M^4.$ Thus  $D_{4,4}^4\neq 0. $
\bigskip

{\bf 4. $D_{4,4}^6\neq 0$:}
By Lemma \ref{vectors}, we have
$b=e^{-\e_2-\e_3+\e_{17}+e_{18}}\in M^4. $
One immediately see that
$$b_0 v^2=C_2 e^{\e_1-\e_3+\e_{17}+\e_{18}-\e_{19}-\e_{20}}$$
 for some $C_2\in \{\pm 1\}.$
We need to show the projection of $b_0 v^2$ to $M^6$ is nonzero.
It suffices to show $M^0\cdot b_0 v^2=\<u_nb_0 v^2|u\in M^0,n\in\Z\>$ has a nonzero projection to $M^6.$
Applying the operator $e_0e_0e_0$ to  $b_0 v^2,$ we get
\begin{equation*}
(e_0)(e_0)(e_0) b_0 v^2=C_3 e^{\e_1+\e_7+\e_8-\e_{13}-\e_{19}-\e_{20}} \in M^0\cdot b_0 v^2
\end{equation*}
for some  $C_3 \ne 0.$  Since $e^{\e_1+\e_7+\e_8-\e_{13}-\e_{19}-\e_{20}}\in M^6$ (Lemma \ref{vectors}),
 the projection of $M^0\cdot b_0 v^2 $ to $M^6$ is nonzero, i.e. $D_{4,4}^6\neq 0.$

\end{document}